\documentclass[11pt]{article}

\usepackage[T1]{fontenc}
\usepackage{lmodern}
\usepackage[margin=1in]{geometry}
\usepackage{amsmath,amssymb,amsthm,mathtools}
\usepackage{booktabs,array,enumitem,microtype,tikz,young}
\usepackage[hidelinks]{hyperref}

\newtheorem{theorem}{Theorem}[section]
\newtheorem{lemma}[theorem]{Lemma}
\newtheorem{proposition}[theorem]{Proposition}
\newtheorem{corollary}[theorem]{Corollary}
\theoremstyle{definition}
\newtheorem{definition}[theorem]{Definition}
\newtheorem{example}[theorem]{Example}
\theoremstyle{remark}
\newtheorem{remark}[theorem]{Remark}

\newcommand{\Asc}{\operatorname{Asc}}
\newcommand{\Des}{\operatorname{Des}}
\newcommand{\Occ}{\operatorname{Occ}}
\newcommand{\NF}{\operatorname{NF}}
\newcommand{\sgn}{\operatorname{sgn}}
\newcommand{\cA}{\mathcal A}
\newcommand{\cD}{\mathcal D}
\newcommand{\cO}{\mathcal O}
\newcommand{\cX}{\mathcal X}
\newcommand{\mgraph}{\mathbf m}
\newcommand{\ngraph}{\mathbf n}

\newcommand{\PRBS}{P_{\mathsf{rB}}}
\newcommand{\PCBS}{P_{\mathsf{cB}}}
\newcommand{\PRBSp}{P_{\mathsf{rB}}^{+}}
\newcommand{\PCBSm}{P_{\mathsf{cB}}^{-}}
\newcommand{\fromRSK}{\xleftarrow{\,\mathsf{RSK}\,}}
\newcommand{\fromRBS}{\xleftarrow{\,\mathsf{rB}\,}}
\newcommand{\fromCBS}{\xleftarrow{\,\mathsf{cB}\,}}
\newcommand{\rowtableau}[1]{\vcenter{\hbox{\begin{Young}#1\end{Young}}}}

\title{Molecules of Type-$B$ and Odd Type-$D$ Gelfand
$\mgraph$- and $\ngraph$-Graphs}
\author{Yifeng Zhang
\\
School of Mathematical Sciences \\
South China Normal University \\
{\tt calvinz314159@gmail.com}}
\date{}

\begin{document}
\maketitle

\begin{abstract}
Gelfand $W$-graphs provide multiplicity-free canonical models for classical
Weyl groups, and their molecules record the connectivity generated by
bidirected canonical-basis edges.  Restricted Beissinger insertion and dual
equivalence give a natural tableau language for these edges, while the wall
and fork generators require additional local criteria.  For both type-$B$
Gelfand $\mgraph$- and $\ngraph$-graphs, we prove that the molecules are exactly
the connected components of explicit occurrence--wall graphs constructed
from split modified Beissinger insertion.  In odd type $D$, we likewise
classify the $\mgraph$-molecules by a completion-filtered occurrence--fork
graph and obtain the $\ngraph$-classification through the parity-sensitive
type-$D$ duality.
\end{abstract}

\section{Introduction}

A Gelfand model is a representation containing every irreducible
representation with multiplicity one.  Perfect models provide a uniform
Coxeter-theoretic framework for such representations in the classical
Weyl groups \cite{MZPerfect}.  Their Hecke deformations are naturally
connected with Kazhdan--Lusztig $W$-graphs and canonical bases
\cite{KazhdanLusztig,MarbergBar}.  In particular, Marberg and Zhang
constructed two companion Gelfand $W$-graphs for the classical Weyl
groups, denoted here by the bold upright symbols $\mgraph$ and $\ngraph$,
and realized their vertices through model triples and fixed-point-free
completions \cite{MZGelfand}.  In type $B$, the vertices of both graphs
may be represented concretely by signed involutions.  Their canonical
bases are related by sign opposition and Hecke contragredience.  This
duality reverses arrows and therefore preserves the bidirected relation,
but it changes the wall convention; keeping track of that change is
essential for an explicit classification.

Kazhdan--Lusztig $W$-graphs for type $B$ arise from the Hecke algebra of
the hyperoctahedral group.  In the unequal-parameter setting their form
depends on a weight function, and the transition between the equal,
intermediate, and asymptotic regimes is a central feature of the theory
\cite{LusztigUnequal,BonnafeIancu}.  Tableau models use domino insertion
and generalized Robinson--Schensted correspondences.  In the geometric
approach of Garfinkle, moving through open cycles of a domino tableau
controls the relevant tableau equivalences \cite{Garfinkle}; related
equivalence classes for type $B$ were developed further in
\cite{Pietraho}.  In the asymptotic parameter regime, Bonnaf\'e and Iancu
gave a generalized Robinson--Schensted description
\cite{BonnafeIancu}.

The subgroup $D_n\subset B_n$ consists of the even signed permutations
and has the forked Coxeter generating set
$s_{-1},s_1,s_2,\ldots,s_{n-1}$.  Type $D$ is not a formal restriction of
the type-$B$ graph: the extra generator $s_{-1}$ is not the type-$B$ wall
reflection, and its strict and weak ascent conditions are computed from a
fixed-point-free completion.  Moreover, perfect models exist in type $D$
only in odd rank \cite{MZPerfect}.  In those ranks the associated canonical
graphs are Gelfand $W$-graphs.  The two type-$D$ graphs are again dual, but
the duality includes the diagram automorphism $s_{-1}\leftrightarrow s_1$
on parity-selected summands \cite{MZGelfand}.  These facts explain both why
the type-$B$ insertion theory remains useful and why a separate fork
analysis is necessary.

The special role of the fork also has a precedent in the ordinary
Kazhdan--Lusztig cell theory of type $D$.  Garfinkle and Vogan developed
operators adapted to branched Dynkin diagrams \cite{GarfinkleVogan}, while
the domino-tableau approach of Garfinkle and McGovern describes left cells
in the classical Weyl groups and records the additional type-$D$
modification \cite{Garfinkle,McGovernDomino}.  Those results concern cells in
the regular Kazhdan--Lusztig graph.  The graphs studied here are instead the
canonical $W$-graphs of multiplicity-free Gelfand models.  Thus the
occurrence--fork theorem below is not a reformulation of the classical
domino-cell classification: it determines exactly which strict
$s_{-1}$-conjugations remain bidirected after the positive-color tail has
been contracted.

A \emph{molecule} of a $W$-graph is a connected component after only the
bidirected edges are retained.  This local viewpoint is also fundamental
in the theory of admissible $W$-graphs \cite{Stembridge}.  For the usual
type-$B$ Kazhdan--Lusztig graphs, the tableau descriptions above organize
molecules by local dual-equivalence and open-cycle moves; the precise
form depends on the parameter regime.  The Gelfand $\mgraph$- and
$\ngraph$-graphs are different canonical graphs on multiplicity-free
models, so their molecules require a separate analysis rather than a
direct transfer of the domino classification.

Insertion algorithms supply the other ingredient.  Beissinger gave a
tableau construction for involutions \cite{Beissinger}; Haiman's dual
equivalence organizes standard tableaux by elementary three-letter moves
\cite{Haiman}, with a graph-theoretic development in \cite{Assaf}.  The
restricted row and column variants needed for Gelfand graphs were
developed in \cite{MZInsertion}.  The present paper combines these local
insertion moves with the type-$B$ wall reflection.

The type-$B$ classification proceeds in two stages.  First split a signed
involution into its positive and negative carriers, apply restricted row
and column Beissinger insertion, and contract the components generated by
the resulting local tableau moves.  The remaining bidirected edges are
wall edges and are detected by the two closed conditions proved in
Section~\ref{sec:wall}.  For the companion $\ngraph$-graph, sign opposition
exchanges the two insertion roles and reverses the wall convention, but
leaves the same two-stage structure.  The two main type-$B$ statements can
therefore be recorded together as follows.

\begingroup
\renewcommand{\thetheorem}{\ref{thm:classification} and
\ref{thm:n-classification}}
\begin{theorem}[Type-$B$ molecule classifications]
Let $y,z$ be signed involutions in the same fixed two-cycle block.
\begin{enumerate}[label=\textup{(\alph*)}]
\item The vertices $y,z$ belong to the same type-$B$ Gelfand
$\mgraph$-molecule if and only if their modified Beissinger occurrences
lie in the same connected component of the $\mgraph$ occurrence--wall graph.
\item The vertices $y,z$ belong to the same type-$B$ Gelfand
$\ngraph$-molecule if and only if their split $\ngraph$-occurrences lie in
the same connected component of the $\ngraph$ occurrence--wall graph.
\end{enumerate}
\end{theorem}
\endgroup

For type $D$ the result takes the following completion-based form.

\begingroup
\renewcommand{\thetheorem}{\ref{thm:D-classification}}
\begin{theorem}[Type-$D$ molecule classifications]
Let $n$ be odd.  Two vertices in the same type-$D_n$ model block belong to
the same Gelfand $\mgraph$-molecule if and only if their
completion-filtered type-$D$ tail occurrences lie in the same connected
component of the
type-$D$ occurrence--fork graph.  The parity-sensitive type-$D$ duality
transports these components bijectively to the Gelfand
$\ngraph$-molecules.
\end{theorem}
\endgroup


Section~\ref{sec:carriers} fixes the signed-carrier and completion
conventions and recalls the canonical bidirected-edge theorem.
Section~\ref{sec:wall} proves the closed wall criterion.
Section~\ref{sec:insertion} defines split modified Beissinger insertion,
works out representative examples, and proves the positive local-carrier
theorem.  Section~\ref{sec:classification} constructs the
occurrence--wall graph, proves the classification, and records the
resulting algorithm.  Section~\ref{sec:n-case} develops the companion
$\ngraph$-classification from sign-opposition duality and states it in
the same occurrence--wall form.  Section~\ref{sec:type-d} restricts the
carrier model to type $D$, proves the fork criterion by a general
commuting-label lemma, constructs the occurrence--fork graph, and obtains
the $\ngraph$-classification from the modified type-$D$ duality.

\section{Preliminaries}
\label{sec:carriers}

\subsection{Signed permutations and the Coxeter groups of types
\texorpdfstring{$B$ and $D$}{B and D}}

Let $B_n$ be the group of signed permutations of $[n]=\{1,\ldots,n\}$.
Thus an element is a bijection of $\{\pm1,\ldots,\pm n\}$ satisfying
$w(-i)=-w(i)$.  Its standard Coxeter generators are $s_0,s_1,\ldots,s_{n-1}$:
$s_0$ changes the sign of $1$, and $s_i$ for $i>0$ interchanges $i,i+1$
and simultaneously interchanges $-i,-i-1$.  An edge associated with
$s_i$ has \emph{color} $i$; color $0$ is the \emph{wall color}, while
the colors $1,\ldots,n-1$ are the \emph{positive colors}.

For $n\geq3$, the type-$D$ Coxeter group is the subgroup
\[
 D_n=\{w\in B_n:|\{i\in[n]:w(i)<0\}|\text{ is even}\}.
\]
Its simple generators are
\[
 s_{-1}=(1,-2)(-1,2),\qquad s_1,s_2,\ldots,s_{n-1}.
\]
The two fork generators $s_{-1}$ and $s_1$ commute, and the only simple
generator that fails to commute with $s_{-1}$ is $s_2$.  We refer to
$s_{-1}$ as the \emph{extra fork color}; the colors
$s_1,\ldots,s_{n-1}$ form the type-$A$ tail.

An involution $z\in B_n$ is written as a signed matching.  Its components
are
\[
 A_i^+ : z(i)=i,\qquad A_i^- : z(i)=-i,
\]
and, for $i<j$,
\[
 Z_{ij}^+ : z(i)=j,\ z(j)=i,
 \qquad
 Z_{ij}^- : z(i)=-j,\ z(j)=-i.
\]
The sign of every letter is the sign of its component.  Let
$\cX_{n,k}$ be the signed involutions with exactly $k$ pair components.
Every carrier move defined below preserves $k$, so $\cX_{n,k}$ is a block
for all the graphs considered below.

For $i>0$, let $\tau_i z=s_i z s_i$ when this conjugate differs from $z$.
The wall move $\tau_0$ is defined by changing the sign of the component
containing $1$.  We refer to the fixed points and pairs collectively as
\emph{carriers}.

\subsection{The type-\texorpdfstring{$B$}{B} Gelfand
\texorpdfstring{$\mgraph$}{m}- and \texorpdfstring{$\ngraph$}{n}-graphs}

Let $\cA=\mathbb Z[v,v^{-1}]$, and let $\mathcal H(B_n)$ be the
equal-parameter Iwahori--Hecke algebra with generators
$H_0,H_1,\ldots,H_{n-1}$.  They satisfy the braid relations of $B_n$ and
the quadratic relations
\[
                      (H_i-v)(H_i+v^{-1})=0.
\]
The two type-$B$ Gelfand modules of \cite{MZGelfand} are
$\mathcal H(B_n)$-modules with standard bases indexed, block by block, by
the same sets $\cX_{n,k}$.  Their bar involutions and triangular orders
determine lower canonical bases
\[
 \{C_z^{\mgraph}:z\in\cX_{n,k}\},
 \qquad
 \{C_z^{\ngraph}:z\in\cX_{n,k}\}.
\]
For $\epsilon\in\{\mgraph,\ngraph\}$, the Gelfand
$\epsilon$-graph has an arrow $z\to y$ whenever some $H_iC_z^\epsilon$
has a nonzero off-diagonal $C_y^\epsilon$-coefficient.  An
$\epsilon$-\emph{molecule} is a connected component of the undirected
graph whose edges are precisely the pairs for which both directed arrows
occur.

Write $D_{\mgraph}(z)$ and $D_{\ngraph}(z)$ for the descent-label sets of
the two graphs.  If $-z$ denotes the signed involution obtained by
reversing the sign of every carrier, sign-opposition duality gives
\begin{equation}
 D_{\ngraph}(z)=[0,n-1]\setminus D_{\mgraph}(-z).
 \label{eq:mn-label-duality}
\end{equation}
It also reverses every directed off-diagonal coefficient; the precise
statement and proof are given in Lemma~\ref{lem:mn-duality}.  Sections
\ref{sec:wall}--\ref{sec:classification} treat the $\mgraph$-graph, and
there we abbreviate $C_z^{\mgraph}$ to $C_z$ and
$D_{\mgraph}(z)$ to $D(z)$.  Section~\ref{sec:n-case} translates the
entire construction to the $\ngraph$-graph.

\subsection{Signed involutions and fixed-point-free completion}

The completion convention is needed to determine strict and weak
ascents.  List the negative fixed points of $z$ as
$c_1>\cdots>c_p$ and the positive fixed points as
$e_1<\cdots<e_q$.  Define the fixed-point-free signed completion
$v=\iota(z)$ on $\{\pm1,\ldots,\pm2n\}$ by $v(-a)=-v(a)$ and
\begin{align}
 v(c_j)&=-(n+j),&v(n+j)&=-c_j &&(1\leq j\leq p),\label{eq:neg-completion}\\
 v(e_j)&=n+p+j,&v(n+p+j)&=e_j &&(1\leq j\leq q),\label{eq:pos-completion}
\end{align}
together with $v(i)=z(i)$ at every nonfixed $i\leq n$.  Pair the unused
positive auxiliary letters consecutively and positively.  This last
choice does not affect any comparison at a node $0,1,\ldots,n-1$.

A node $i$ is called \emph{strict} at $z$ when the carrier move
$\tau_i z$ is nontrivial; otherwise it is \emph{weak}.  The canonical
model orders the two endpoints of every strict move as its lower and upper
endpoint.  We only use this published orientation when invoking the
strict-cover theorem; all classification graphs below are undirected.

Write $\Asc(z)$ for the ascent label of the type-$B$ Gelfand $\mgraph$-graph
and put
\[
                         D(z)=[0,n-1]\setminus\Asc(z).
\]
The following formula makes the convention intrinsic to the signed
matching.

\begin{lemma}[Completion descents and strict conjugates]
\label{lem:completion}
For $v=\iota(z)$ one has
\begin{align}
 D(z)={}&\{0:v(1)<0\}\notag\\
 &\cup\{i\in[n-1]:v(i)=\pm(i+1)\}\notag\\
 &\cup\{i\in[n-1]:i,i+1\text{ are not same-sign fixed points and }
                                  v(i)>v(i+1)\}.\label{eq:descent-formula}
\end{align}
For $i>0$, the conjugate $\tau_i z$ is nontrivial and strict exactly when
$i,i+1$ are neither the endpoints of one pair nor same-sign fixed points.
The carrier $\tau_0z$ is obtained by reversing the sign of the component
containing $1$, and is always strict.
\end{lemma}

\begin{proof}
The equality $v(i)=\pm(i+1)$ is precisely the weak-descent case.  External
completion values occur at original letters only at fixed points.  The
decreasing order in \eqref{eq:neg-completion} and the increasing order in
\eqref{eq:pos-completion} make two consecutive external values increasing
precisely for two same-sign fixed points; this is the weak-ascent case.
All remaining nodes are strict, and their direction is determined by the
comparison of $v(i)$ and $v(i+1)$.  The node-zero formula is the defining
sign comparison.

Interchanging $i,i+1$ fixes the matching exactly when they are the
endpoints of one pair or are same-sign fixed points.  Otherwise it
relabels distinct carriers or distinct endpoints and is strict.  The
$s_0$ assertion is immediate when $1$ lies in a pair.  If $1$ is fixed,
then it is the last member of the decreasing negative list when negative,
and the first member of the increasing positive list when positive.
Changing its sign therefore
uses the same auxiliary partner and leaves every other partner unchanged.
Hence $s_0\iota(z)s_0=\iota(\tau_0z)$ in the fixed-point case as well.
\end{proof}

\subsection{Type-\texorpdfstring{$D$}{D} carrier blocks and graph labels}
\label{subsec:type-D-prelim}

A carrier is negative when its superscript is $-$.  Set
\begin{equation}
 \cX^D_{n,k}=\{z\in\cX_{n,k}:
       \text{the number of negative carriers of $z$ is even}\}.
 \label{eq:type-D-carriers}
\end{equation}
The completion $v=\iota(z)$ defined in
\eqref{eq:neg-completion}--\eqref{eq:pos-completion} maps
$\cX^D_{n,k}$ bijectively onto the corresponding type-$D$ model block;
this is the signed-carrier form of \cite[(2.3)]{MZGelfand}.

For $i\geq1$, membership of $s_i$ in the type-$D$ descent label
$D_D(z)$ is given by the positive-node part of
\eqref{eq:descent-formula}.  At the extra fork node,
\cite[Proposition~2.5(c)]{MZGelfand} gives four mutually exclusive cases:
\begin{align*}
 s_{-1}&\text{ is a weak descent} &&\Longleftrightarrow v(1)=\pm2,\\
 s_{-1}&\text{ is a weak ascent}  &&\Longleftrightarrow
 \begin{cases}
 v(1)<-n<n<v(2),\quad\text{or}\\
 v(2)<-n<n<v(1),
 \end{cases}\\
 s_{-1}&\text{ is a strict descent}&&\Longleftrightarrow
 \text{neither weak case holds and }-v(2)>v(1),\\
 s_{-1}&\text{ is a strict ascent}&&\Longleftrightarrow
 \text{neither weak case holds and }-v(2)<v(1).
\end{align*}
The descent label $D_D(z)$ is the union of the weak and strict descents.
When $s_{-1}vs_{-1}\ne v$, define the extra fork move by
\begin{equation}
 \tau^D_{-1}z=\iota^{-1}(s_{-1}\iota(z)s_{-1}).
 \label{eq:type-D-fork-toggle}
\end{equation}
This is deliberately a completion-level definition.  The auxiliary
partners of fixed carriers depend on all fixed carriers of the same sign,
so equality of the first three completion values is not a local
three-carrier statement.

Let $\Gamma^{\mgraph}_{D,n,k}$ and $\Gamma^{\ngraph}_{D,n,k}$ denote the
two canonical type-$D$ graphs on this block.  They are Gelfand $W$-graphs
when $n$ is odd; no perfect model exists in even type-$D$ rank
\cite{MZPerfect}.  Let $\diamond$ be the diagram automorphism interchanging
$s_{-1}$ and $s_1$ and fixing the other simple generators.  We write
$\delta_n^D$ for the parity-sensitive involution denoted $\iota_n^D$ in
\cite[(3.26)]{MZGelfand}; on each model block it is sign opposition,
possibly followed by $\diamond$.  The type-$D$ duality theorem
\cite[Theorem~1.13]{MZGelfand} identifies
$\Gamma^{\mgraph}_{D,n,k}$ with the dual of a $\diamond$-modified
$\ngraph$-graph, and the modified graph has the same underlying directed
graph as $\Gamma^{\ngraph}_{D,n,k}$.

\subsection{The canonical bidirected-edge theorem}

The completion construction above is the signed-matching form of the
model-triple parametrization in \cite{MZGelfand}.  More precisely,
\cite[Lemmas~3.27--3.28 and Theorem~1.10]{MZGelfand} identify
$z\mapsto\iota(z)$ with the relevant fixed-$k$ model-triple block and
transport the height order, the four strict/weak ascent and descent
types, and conjugation by each simple generator.  Thus the signed-carrier
labels and the canonical $W$-graph labels used here are the same labels,
rather than merely sets of the same cardinality.

We use the following specialization of the strict-cover theorem from
\cite[Corollary~3.14]{MZGelfand}.  It is restated here to isolate exactly
what enters the classification.

\begin{theorem}[Canonical bidirected-edge theorem]
\label{thm:all-edges}
In either a type-$B$ block $\cX_{n,k}$ or a type-$D$ block
$\cX^D_{n,k}$, two vertices $z,z'$ form a bidirected edge in the
corresponding $\mgraph$-graph if and only if $z'$ is a nontrivial strict
simple conjugate of $z$ and
\begin{equation}
 \text{their two descent labels are incomparable}.
\label{eq:strict-edge}
\end{equation}
Every such edge has coefficient one.  In particular, there are no
additional nonlocal bidirected canonical-basis edges.
\end{theorem}

\begin{proof}
Write $E(u)$ for the appropriate type-$B$ or type-$D$ descent label and
$A(u)$ for its ascent complement.  
Orient the strict cover from its lower endpoint $z$ to its upper endpoint
$z'$.  The
bidirected-edge theorem of \cite[Corollary~3.14]{MZGelfand} says that the
reverse coefficient is nonzero exactly when
\[
                         E(z)\cap A(z')\ne\varnothing.
\]
The conjugating simple generator is an ascent of $z$ and a descent of
$z'$; the coefficient-one
assertion is the strict
action formula in \cite[Lemma~3.11 and Theorem~3.13]{MZGelfand}.
Therefore the displayed intersection is
nonempty exactly when each of $E(z),E(z')$ contains an element absent
from the other, which is \eqref{eq:strict-edge}.  Lemma~\ref{lem:completion}
identifies the type-$B$ strict covers with the stated signed-carrier
conjugates.  In type $D$, the corresponding identification and the fork
labels are \cite[Lemmas~3.27--3.28 and Proposition~2.5]{MZGelfand},
together with \eqref{eq:type-D-fork-toggle}.
The rank-one type-$B$ case is immediate.
\end{proof}

\section{The wall reflection}
\label{sec:wall}

The wall toggle is $z\leftrightarrow\tau_0z$, namely the sign change of
the carrier through $1$.  Therefore, we have the following criterion:

\begin{theorem}[Wall edge criterion]
\label{thm:wall}
The wall toggle $z\leftrightarrow\tau_0z$ is bidirected if and only if
one of the following mutually exclusive conditions holds.
\begin{enumerate}[leftmargin=3.2em]
\item[\textup{(WF)}]
The component containing $1$ is a fixed point, $n\geq2$, and $2$ is not
a positive fixed point.
\item[\textup{(WP)}]
The component containing $1$ is $Z_{1r}^{\epsilon}$ and the component
containing $2$ is $Z_{2q}^{\delta}$ with $2<q<r$.  The signs
$\epsilon,\delta$ are arbitrary.
\end{enumerate}
In particular, there is no wall edge for $n=1$ or for $Z_{12}^{\epsilon}$.
\end{theorem}

\begin{proof}
Orient the toggle from a negative component through $1$ to a positive one
and call the two matchings $z^-,z^+$.  Then
$0\in D(z^-)\setminus D(z^+)$.  By Theorem~\ref{thm:all-edges}, the
toggle is bidirected precisely when some nonzero descent is created.

Suppose first that $1$ is fixed.  In $z^-$ its value is the smallest
negative external value, whereas in $z^+$ its value is larger than every
internal value.  The auxiliary partner of $1$ and all other auxiliary
partners are unchanged, as observed in Lemma~\ref{lem:completion}.  Thus
only node $1$ can gain a descent.  If $2$ is a positive fixed point, node
$1$ is an ascent on both sides.  If $2$ is a negative fixed point, it
changes from a weak ascent to a strict descent.  If $2$ is nonfixed, it
changes from a strict ascent to a strict descent.  This proves
\textup{(WF)}, including the condition $n\geq2$.

Now suppose that $1$ is paired with $r$.  Only the four values
\[
 v^-(1)=-r,\quad v^-(r)=-1,
 \qquad
 v^+(1)=r,\quad v^+(r)=1
\]
change.  If $r=2$, node $1$ is a weak descent before and after the toggle.
At node $2$, when it exists, a new descent would require the nonzero
integer $v(3)$ to lie strictly between $-1$ and $1$, which is impossible.
The descent sets are therefore comparable.

Assume $r>2$.  At $r-1$, replacing the second comparison value $-1$ by
$1$ can only remove a descent.  At $r$ (absent when $r=n$), replacing the
first comparison value can create a descent only if $v(r+1)$ lies strictly
between $-1$ and $1$, again impossible.  A new descent can therefore occur
only at node $1$, and it occurs exactly when
\[
                              -r<v(2)<r.
\]
A fixed point at $2$ has an external value of absolute value greater than
$n$.  If $2$ is paired with $q>2$, then $v(2)=\pm q$, so the inequality
is equivalent to $q<r$.  This proves \textup{(WP)}, including the
boundary case $r=n$.
\end{proof}

\begin{proposition}[Endpoint tableau form of the wall criterion]
\label{prop:tableau-wall-criterion}
Let $z^-,z^+$ be the two endpoints of a wall toggle, labelled so that the
carrier through $1$ is negative in $z^-$ and positive in $z^+$.  Then
\begin{equation}
 z^-\longleftrightarrow z^+
 \quad\Longleftrightarrow\quad
 \Des\bigl(\Psi(z^+)\bigr)
       \setminus\Des\bigl(\Psi(z^-)\bigr)\ne\varnothing.
 \label{eq:tableau-wall-criterion}
\end{equation}
Equivalently, if
\[
 \Des_B^{\mathrm{tab}}(\Psi(u))=
 \Des(\Psi(u))\cup
 \begin{cases}
  \{0\},&1\text{ lies in the negative-transpose component of }\Psi(u),\\
  \varnothing,&1\text{ lies in the positive component of }\Psi(u),
 \end{cases}
\]
then the wall edge is bidirected exactly when the two augmented tableau
descent sets are incomparable.  This is an endpoint test: it does not
claim that $\Psi(z^+)$ can be constructed from $\Psi(z^-)$ alone.
\end{proposition}

\begin{proof}
Lemma~\ref{lem:tableau-descents} identifies the positive-node parts of
the augmented tableau descent sets with the canonical descent labels,
while $0\in D(z^-)\setminus D(z^+)$.  Thus
\eqref{eq:tableau-wall-criterion} is precisely the incomparability test
of Theorem~\ref{thm:all-edges}.
\end{proof}

\begin{example}[The two wall mechanisms]
The edge
\[
 A_1^-Z_{24}^+Z_{36}^+A_5^-
 \longleftrightarrow
 A_1^+Z_{24}^+Z_{36}^+A_5^-
\]
is of type \textup{(WF)}, since $1$ is fixed and $2$ is not a positive
fixed point.  Its wall toggle is visible in the following arc diagram;
the sign above an isolated node is the sign of its fixed carrier.
\begin{center}
\begin{tikzpicture}[x=.55cm,y=.55cm,
  v/.style={circle,draw,fill=white,inner sep=0pt,minimum size=3.7mm},
  every path/.style={thick}]
  \foreach \x/\lab in {0/1,1/2,2/3,3/4,4/5,5/6}
    \node[v] (a\lab) at (\x,0) {\scriptsize\lab};
  \draw (a2) to[bend left=48] node[above,font=\scriptsize] {$+$} (a4);
  \draw (a3) to[bend left=55] node[above,font=\scriptsize] {$+$} (a6);
  \node[font=\scriptsize] at (0,.58) {$-$};
  \node[font=\scriptsize] at (4,.58) {$-$};
\end{tikzpicture}
\quad $\longleftrightarrow$ \quad
\begin{tikzpicture}[x=.55cm,y=.55cm,
  v/.style={circle,draw,fill=white,inner sep=0pt,minimum size=3.7mm},
  every path/.style={thick}]
  \foreach \x/\lab in {0/1,1/2,2/3,3/4,4/5,5/6}
    \node[v] (b\lab) at (\x,0) {\scriptsize\lab};
  \draw (b2) to[bend left=48] node[above,font=\scriptsize] {$+$} (b4);
  \draw (b3) to[bend left=55] node[above,font=\scriptsize] {$+$} (b6);
  \node[font=\scriptsize] at (0,.58) {$+$};
  \node[font=\scriptsize] at (4,.58) {$-$};
\end{tikzpicture}
\end{center}
The same two mechanisms can now be read directly from the endpoint
tableaux.  The smallest fixed-carrier example is
\[
\begin{array}{c@{\quad}c@{\quad}c}
z&\Psi(z)&\Des(\Psi(z))\\ \hline
A_1^-A_2^-&
  \varnothing\sqcup\rowtableau{1&2\cr}&\varnothing\\[2pt]
A_1^+A_2^-&
  \rowtableau{1\cr}\sqcup\rowtableau{2\cr}&\{1\}.
\end{array}
\]
The negative-transpose component is displayed second but is read before
the positive component.  Hence the wall toggle creates the tableau
descent $1$, as required by \eqref{eq:tableau-wall-criterion}.

A wall edge of type \textup{(WP)} has the schematic form
\[
 Z_{1r}^{\epsilon}Z_{2q}^{\delta}
 \longleftrightarrow
 Z_{1r}^{-\epsilon}Z_{2q}^{\delta},
 \qquad 2<q<r,
\]
with every other carrier unchanged.  Diagrammatically, the arc through
$2$ is properly nested inside the arc through $1$.
\begin{center}
\begin{tikzpicture}[x=1cm,y=.58cm,
  v/.style={circle,draw,fill=white,inner sep=0pt,minimum size=4.2mm}]
  \node[v] (p1) at (0,0) {\scriptsize$1$};
  \node[v] (p2) at (1,0) {\scriptsize$2$};
  \node[v] (pq) at (3,0) {\scriptsize$q$};
  \node[v] (pr) at (4.5,0) {\scriptsize$r$};
  \draw[very thick] (p1) to[bend left=55]
    node[above,font=\scriptsize] {$\epsilon\mapsto-\epsilon$} (pr);
  \draw[thick] (p2) to[bend left=42]
    node[above,font=\scriptsize] {$\delta$} (pq);
  \node[font=\scriptsize] at (2.25,-.62) {$1<2<q<r$};
\end{tikzpicture}
\end{center}
For instance, $Z_{14}^-Z_{23}^+\leftrightarrow
Z_{14}^+Z_{23}^+$ gives the following illustrated descent test:
\[
\begin{array}{c@{\qquad}c@{\qquad}c}
z&\Psi(z)&\Des(\Psi(z))\\ \hline
Z_{14}^-Z_{23}^+&
 \rowtableau{2\cr3\cr}\sqcup\rowtableau{1\cr4\cr}&\{2,3\}\\[2pt]
Z_{14}^+Z_{23}^+&
 \rowtableau{1\cr2\cr3\cr4\cr}\sqcup\varnothing&\{1,2,3\}.
\end{array}
\]
Thus $1$ belongs to the upper endpoint's tableau descent set but not to
the lower endpoint's, so the nested-arc wall toggle is bidirected.
\end{example}

\section{Split modified Beissinger insertion}
\label{sec:insertion}

\subsection{Restricted row and column insertion}

A \emph{partially standard tableau} is a Young tableau with distinct
positive-integer entries, increasing along rows and columns.  If $a$ is
not an entry of such a tableau $T$, then $T\fromRSK a$ denotes ordinary
row insertion.  Suppose $a\leq b$ are absent from $T$.  If $a=b$, set
\[
                 T\fromRBS(a,a)=T\fromCBS(a,a)=T\fromRSK a.
\]
If $a<b$, let $(r,c)$ be the box of $T\fromRSK a$ that is not in $T$.

Following the notation of the row Beissinger correspondence, define
\[
 T\fromRBS(a,b)
\]
by appending $b$ to row $r+1$ of $T\fromRSK a$.  Define
\[
 T\fromCBS(a,b)
\]
by appending $b$ to the bottom of column $c+1$ of $T\fromRSK a$.
These are the row and column Beissinger pair insertions
\cite{Beissinger,MZInsertion}.

For a fixed-point-free involution on an ordered alphabet, its \emph{full
row} (respectively, \emph{full column}) Beissinger tableau is obtained by
applying $\mathsf{rB}$ (respectively, $\mathsf{cB}$) to all two-cycles in
increasing order of their larger endpoints.  We denote the resulting
tableaux by $\PRBS$ and $\PCBS$, respectively.

For later use, we define the \emph{three-letter betweenness operation}.
Suppose a word contains three distinguished letters $a<b<c$.  If $b$
occurs positionally between $a$ and $c$, the operation fixes the word.
Otherwise, it exchanges $b$ with whichever of $a$ and $c$ is farther from
$b$ in the word.  Thus the operation depends only on the relative order of
the three distinguished positions.

\begin{lemma}[Restriction and standardization]
\label{lem:restriction-standardization}
Let $S$ be a finite ordered alphabet.
\begin{enumerate}[label=\textup{(\alph*)}]
\item Order-preserving relabeling of $S$ commutes with row insertion,
row and column Beissinger pair insertion, their inverse operations, and
the three-letter betweenness operation defined above.
\item Complete an involution on $S$ by pairing its fixed points with new
letters larger than every member of $S$.  If the auxiliary right endpoints
follow the fixed points in increasing order, then restricting the full row
Beissinger tableau back to $S$ amounts to first processing the nontrivial
pairs and then row-inserting the fixed points in increasing order.  If the
auxiliary right endpoints follow the fixed points in decreasing order,
then the analogous restricted column construction row-inserts the fixed
points in decreasing order.
\end{enumerate}
Both restricted constructions are bijective on their respective
fixed-point blocks.
\end{lemma}

\begin{proof}
Every comparison in a row-bumping path depends only on the relative order
of the letters.  The row containing the new box, and hence the row or
column to which the second member of a Beissinger pair is appended, is
therefore unchanged by order-preserving relabeling.  Reversing these steps
proves the inverse statement.  The three-letter betweenness operation uses
only the relative order of its entries in the reading word.  This proves
(a).

In the increasing-order completion, all auxiliary right endpoints are larger than
the original alphabet, so their pairs are processed after the original
two-cycles.  For a completed fixed point $c$, pair insertion first
row-inserts $c$ and then appends its auxiliary partner.  Deleting that
partner leaves exactly the row insertion of $c$.  The order of the
auxiliary right endpoints is the increasing order of the fixed points.
This gives the row assertion.  The decreasing-order completion orders the
auxiliary partners oppositely; after deleting each appended auxiliary box,
the column construction leaves row insertion of the fixed points in
decreasing order.  These are precisely the completion-restriction
identities of \cite[Lemmas~3.13 and~3.16]{MZInsertion}.

The inverse restricted constructions and their fixed-point parity
conditions are proved in
\cite[Lemmas~3.13, 3.16 and Theorems~3.14, 3.17]{MZInsertion}.  They
recover a unique completed involution and hence the original involution.
This proves the asserted bijectivity.
\end{proof}

Split the carriers of $z$ by sign.  In either sign submatching, process
the nontrivial pairs $(a,b)$ in increasing order of $b$.  Use
$\mathsf{rB}$ on the
positive submatching and then row-insert its fixed points in increasing
order.  Use $\mathsf{cB}$ on the negative submatching and then row-insert
its fixed points in decreasing order.  Denote the resulting partially
standard tableaux by $\PRBSp(z)$ and $\PCBSm(z)$.  The superscripts record
the sign submatching; they are not powers.  Their alphabets are disjoint
and cover $[n]$.  Applying Lemma~\ref{lem:restriction-standardization} after
order-preserving standardization of the two sign alphabets shows that
\begin{equation}
                      z\longmapsto(\PRBSp(z),\PCBSm(z))
\label{eq:split-insertion}
\end{equation}
is injective and has an explicit inverse on its image.  We call a pair of
partial tableaux \emph{admissible} if it belongs to this image.
The odd-column and odd-row parities recover the positive and negative
fixed-point numbers.

\subsection{The ordered disconnected tableau}

\begin{definition}[Corrected split tableau]
\label{def:Psi}
Define
\[
       \Psi(z)=\PRBSp(z)\sqcup\bigl(\PCBSm(z)\bigr)^{\mathsf t}.
\]
This is an ordered disconnected tableau.  The two diagrams have no
cross-component Young-poset relations.  Each component is read from left
to right along rows and from its bottom row to its top row, and every entry
of the negative-transpose component is declared to precede every entry of
the positive component.
\end{definition}

\begin{example}[A complete split modified Beissinger insertion]
\label{ex:split-insertion}
Consider the signed involution in $B_8$
\[
 z=Z_{14}^+A_2^+A_3^+Z_{58}^-A_6^-A_7^-.
\]
The positive submatching has the pair $(1,4)$ and the fixed points $2,3$.
First apply $\mathsf{rB}$ to $(1,4)$, which gives a column with entries
$1,4$; then
row-insert $2$ and $3$ in increasing order.  Thus
\[
 \PRBSp(z)
 =\emptyset\fromRBS(1,4)\fromRBS(2,2)\fromRBS(3,3)
 =\rowtableau{1&2&3\cr4\cr}.
\]
The negative submatching has the pair $(5,8)$ and the fixed points $6,7$.
The cB-insertion of $(5,8)$ gives the row $5,8$.  Row-inserting the fixed
points in decreasing order, first $7$ and then $6$, gives
\[
 \PCBSm(z)
 =\emptyset\fromCBS(5,8)\fromCBS(7,7)\fromCBS(6,6)
 =\rowtableau{5&6\cr7\cr8\cr},
 \qquad
 \bigl(\PCBSm(z)\bigr)^{\mathsf t}=\rowtableau{5&7&8\cr6\cr}.
\]
Consequently
\[
 \Psi(z)=
 \rowtableau{1&2&3\cr4\cr}
 \sqcup
 \rowtableau{5&7&8\cr6\cr},
\]
where the second displayed component is declared to precede the first in
the ordered row-reading word.  This example also exhibits the two
different completion conventions: positive fixed points are inserted in
increasing order, while negative fixed points are inserted in decreasing
order before transposition.
\end{example}

\begin{example}[Why the negative transpose is necessary]
\label{ex:negative-transpose}
The transpose in Definition~\ref{def:Psi} is forced.  The vertices
\[
 z=Z_{12}^-A_3^+,
 \qquad
 \tau_2z=Z_{13}^-A_2^+
\]
form a bidirected edge.  Without transposition the first local reading
word is $1,2,3$, on which the relevant three-letter move is trivial.  With
the negative tableau transposed, the word is $2,1,3$, and the move gives
the required interchange.
\end{example}

For $2\leq r\leq n-1$, let $d_r$ be Haiman's elementary
dual-equivalence involution on $r-1,r,r+1$ \cite{Haiman}: it fixes the
tableau when $r$ lies between the other two entries in row-reading order,
and otherwise interchanges $r$ with the farther of the two.  More
concretely, after deleting all other letters from the row-reading word, the
orders
\[
 (r-1)\,r\,(r+1)\quad\hbox{and}\quad(r+1)\,r\,(r-1)
\]
are fixed, while the other four orders are paired by
\[
 \begin{array}{c@{\quad\longleftrightarrow\quad}c}
 r\,(r-1)\,(r+1)&(r+1)\,(r-1)\,r,\\
 r\,(r+1)\,(r-1)&(r-1)\,(r+1)\,r.
 \end{array}
\]
The entries not shown retain their positions.  For example,
$d_2(2,1,3)=3,1,2$, whereas $d_2(1,2,3)=1,2,3$.
We also need one boundary move.

\begin{definition}[Disconnected two-letter boundary move]
If $1,2$ lie in different components of $\Psi(z)$, let
$d_0^{\blacktriangledown}$ exchange their entries.  If they lie in the same
component, let $d_0^{\blacktriangledown}$ fix the tableau.
\end{definition}

In music, the solid downward-pointing triangle denotes staccato---a
detached articulation---and here it likewise indicates that the two
entries are disconnected.

\begin{lemma}[Tableau descents]
\label{lem:tableau-descents}
For every signed carrier $z$,
\[
                     D(z)\cap[n-1]=\Des(\Psi(z)),
\]
where $i\in\Des(T)$ means that $i+1$ occurs in a lower row than $i$ in
an ordinary component, equivalently that $i+1$ precedes $i$ in the
declared row-reading word.  The latter condition is the definition when
the two entries lie in different components.
\end{lemma}

\begin{proof}
If $i,i+1$ are both positive, this is the restricted row Beissinger
descent formula in the proof of \cite[Theorem~3.15]{MZInsertion}.  If they
are both negative, negation reverses their completion comparison and
tableau transposition reverses the adjacent tableau descent.  The two
reversals turn the restricted column formula
\cite[Theorem~3.18]{MZInsertion} into the asserted equality for
$\bigl(\PCBSm\bigr)^{\mathsf t}$.

If the signs differ, every negative entry precedes every positive entry in
the reading word of $\Psi(z)$.  Thus $i$ is a tableau descent precisely
when $i$ is positive and $i+1$ is negative.  All positive completion
values are greater than all negative completion values, so this is exactly
the strict-descent condition in \eqref{eq:descent-formula}.
\end{proof}

Consequently the augmented tableau descent set introduced in
\eqref{eq:tableau-wall-criterion} satisfies
\[
                     \Des_B^{\mathrm{tab}}(\Psi(z))=D(z).
\]
Thus the tableau wall test is an exact rewriting of the canonical-label
test, not an additional hypothesis on the endpoints.

\begin{lemma}[The rank-two boundary]
\label{lem:boundary}
If $1,2$ have opposite signs, then
\[
                     \Psi(\tau_1z)=d_0^{\blacktriangledown}\Psi(z),
\]
and $z\leftrightarrow\tau_1z$ is bidirected.
\end{lemma}

\begin{proof}
Each sign alphabet contains exactly one of $1,2$.  Replacing its member by
the other is order-preserving because all remaining letters are at least
$3$.  Lemma~\ref{lem:restriction-standardization}(a) therefore shows that
the tableau swap inverse-inserts to $\tau_1z$.

If the signs at $1,2$ are $+,-$, the descent statuses at nodes $0,1$ are
respectively $(0,1)$; after the swap they are $(1,0)$.  The reverse
orientation is identical.  Thus the descent sets are incomparable, and
Theorem~\ref{thm:all-edges} applies.
\end{proof}

\begin{example}[The indispensable rank-two move]
In $\cX_{2,0}$,
\[
                 A_1^+A_2^-\longleftrightarrow A_1^-A_2^+
\]
has descent sets $\{1\}$ and $\{0\}$.  It is exactly the boundary move
$d_0^{\blacktriangledown}$ and shows why ordinary three-letter dual equivalence
cannot by itself classify the molecules.
\begin{center}
\begin{tikzpicture}[x=.85cm,y=.65cm,
  v/.style={circle,draw,fill=white,inner sep=0pt,minimum size=4.2mm}]
  \node[v] (a1) at (0,0) {\scriptsize$1$};
  \node[v] (a2) at (1,0) {\scriptsize$2$};
  \node[font=\scriptsize] at (0,.55) {$+$};
  \node[font=\scriptsize] at (1,.55) {$-$};
\end{tikzpicture}
\quad $\longleftrightarrow$ \quad
\begin{tikzpicture}[x=.85cm,y=.65cm,
  v/.style={circle,draw,fill=white,inner sep=0pt,minimum size=4.2mm}]
  \node[v] (b1) at (0,0) {\scriptsize$1$};
  \node[v] (b2) at (1,0) {\scriptsize$2$};
  \node[font=\scriptsize] at (0,.55) {$-$};
  \node[font=\scriptsize] at (1,.55) {$+$};
\end{tikzpicture}

{\small The two isolated signed nodes are exchanged across the two
components of the disconnected tableau.}
\end{center}
\end{example}

\subsection{All positive-color edges}

For same-sign letters we use the restricted Beissinger three-carrier
theorem.  It is a direct local insertion result.

\begin{lemma}[Same-sign carrier theorem]
\label{lem:same-sign}
Suppose $r-1,r,r+1$ belong to one sign submatching.  On the positive
submatching use restricted row Beissinger insertion; on the negative
submatching use the transpose of restricted column Beissinger insertion.
Then every nontrivial $d_r$ inverse-inserts to the strict conjugate by
$s_{r-1}$ or $s_r$, according as it exchanges $r-1,r$ or $r,r+1$.
Conversely, every bidirected same-sign strict conjugation involving one of
these adjacent pairs is obtained in this way.
\end{lemma}

\begin{proof}
By Lemma~\ref{lem:restriction-standardization}(a), we may standardize each
sign alphabet independently.  The positive assertion is then the
restricted row theorem \cite[Theorem~3.15]{MZInsertion}; the negative
assertion is the restricted column theorem
\cite[Theorem~3.18]{MZInsertion}.  Elementary dual equivalence commutes
with tableau transposition.  Negating a negative carrier reverses its
strict completion comparisons, while transposition complements adjacent
tableau descents; hence the two conventions agree.  The cited theorems
include all fixed/fixed, fixed/pair, and pair/pair three-carrier cases.
\end{proof}

\begin{theorem}[Positive local-carrier theorem]
\label{thm:positive-local}
The nontrivial moves $d_0^{\blacktriangledown},d_2,\ldots,d_{n-1}$ on the
tableaux $\Psi(z)$ correspond exactly to all bidirected carrier edges of
colors $1,\ldots,n-1$.

More precisely, if $d_r$ exchanges $r-1,r$, inverse insertion gives
$\tau_{r-1}z$, while if it exchanges $r,r+1$, inverse insertion gives
$\tau_rz$.  Conversely, if $z\leftrightarrow\tau_i z$ is bidirected,
then either $i=1$, the signs differ, and this is $d_0^{\blacktriangledown}$, or
\[
 \Psi(\tau_i z)=d_i\Psi(z)
 \quad\text{or}\quad
 \Psi(\tau_i z)=d_{i+1}\Psi(z),
\]
where undefined subscripts are omitted.
\end{theorem}

\begin{proof}
Put $a=r-1,b=r,c=r+1$.  If all three signs agree, the assertion is
Lemma~\ref{lem:same-sign}.  A conjugation is nonstrict precisely when the
two exchanged letters form one pair or are same-sign fixed points; all
other fixed/pair cases are included in that lemma.

Suppose both signs occur.  An entry in the negative-transpose component
precedes every entry in the positive component.  Let $x<_{R}y$ denote
row-reading order within the repeated-sign component.  The complete
mixed-sign calculation is
\[
\begin{array}{c|c|c}
 (\sgn a,\sgn b,\sgn c)&\text{repeated-sign order}&
                           d_r\text{ if nontrivial}\\ \hline
 -++&b<_{R}c&\mathrm{id}\\
 -++&c<_{R}b&(a\ b)\\
 +--&b<_{R}c&(a\ b)\\
 +--&c<_{R}b&\mathrm{id}\\
 +-+&a<_{R}c&(b\ c)\\
 +-+&c<_{R}a&(a\ b)\\
 -+-&a<_{R}c&(a\ b)\\
 -+-&c<_{R}a&(b\ c)\\
 ++-&a<_{R}b&(b\ c)\\
 ++-&b<_{R}a&\mathrm{id}\\
 --+&a<_{R}b&\mathrm{id}\\
 --+&b<_{R}a&(b\ c)
\end{array}
\]
A nontrivial move therefore always exchanges two opposite-sign adjacent
letters.  Each sign submatching changes only by an order-preserving
relabeling, so equivariance of insertion and inverse insertion gives the
corresponding $\tau$-move.  A nontrivial elementary dual-equivalence move
has incomparable tableau descent sets; Lemma~\ref{lem:tableau-descents}
and Theorem~\ref{thm:all-edges} show that this carrier edge is bidirected.

Conversely, strict conjugation by $\tau_i$ reverses the descent at node
$i$ and can change no node outside $i-1,i,i+1$, with node $0$ added only
when $i=1$.  For opposite signs, $\Psi(\tau_i z)$ is the literal swap of
$i,i+1$.  The descent at $i-1$ changes exactly when $i-1$ lies between
them in row-reading order, in which case $d_i$ makes the swap; the
analogous statement at $i+1$ gives $d_{i+1}$.  The exceptional node-zero
change is Lemma~\ref{lem:boundary}.

For equal signs, a neighboring opposite-sign entry is extremal in the
ordered reading word, so its cross-sign comparison cannot change.
Incomparability must therefore occur in a same-sign three-carrier
configuration, already covered by Lemma~\ref{lem:same-sign}.  This proves
both directions.
\end{proof}

\begin{definition}[Modified Beissinger occurrence]
\label{def:occurrence}
Let $\cD_{n,k}$ be the graph on the admissible tableaux $\Psi(z)$ with
$z\in\cX_{n,k}$.  Its edges are all nontrivial
$d_0^{\blacktriangledown},d_2,\ldots,d_{n-1}$.  For a graph $G$ and a
vertex $T$, the notation $[T]_G$ means the vertex set of the connected
component of $G$ containing $T$.  Define
\[
                         \Occ(z)=[\Psi(z)]_{\cD_{n,k}}.
\]
We call this connected component a \emph{modified Beissinger occurrence}.
\end{definition}

Theorem~\ref{thm:positive-local} says that occurrences are exactly the
components obtained by deleting the wall color from the full bidirected
canonical graph.  This is a theorem about explicit tableau moves, rather
than a transported definition.

\section{The occurrence--wall classification}
\label{sec:classification}

\begin{definition}[Occurrence--wall graph]
\label{def:wall-graph}
Let $\cO_{n,k}$ have as vertices the modified Beissinger occurrences in
$\cD_{n,k}$.  Join two occurrences $L,L'$ if there are wall-toggle
endpoints $z^-,z^+\in\cX_{n,k}$ such that
\[
 \Psi(z^-)\in L,\qquad \Psi(z^+)\in L',\qquad
 z^+=\tau_0z^-,
\]
the carrier through $1$ is negative in $z^-$, and
\begin{equation}
 \Des\bigl(\Psi(z^+)\bigr)
       \setminus\Des\bigl(\Psi(z^-)\bigr)\ne\varnothing.
 \label{eq:wall-graph-tableau-edge}
\end{equation}
Multiple endpoint witnesses give one undirected edge.  By
Proposition~\ref{prop:tableau-wall-criterion}, this tableau definition is
equivalent to the closed carrier conditions \textup{(WF)} and \textup{(WP)}.
\end{definition}

Fix any total order on finite tableau pairs.  Let $\NF_A(z)$ be the least
tableau pair in $\Occ(z)$, order occurrences through these representatives,
and let $\NF_B(z)$ be the least occurrence in the connected component of
$\cO_{n,k}$ containing $\Occ(z)$.  The choice of total order only chooses
a name for each component; it does not affect the equivalence relation.

\begin{theorem}[Type-$B$ molecule classification]
\label{thm:classification}
For $y,z\in\cX_{n,k}$, the following are equivalent.
\begin{enumerate}[label=\textup{(\roman*)}]
\item The vertices $y,z$ belong to the same molecule of the type-$B$
Gelfand $\mgraph$-graph.
\item Their modified Beissinger occurrences lie in the same connected
component of $\cO_{n,k}$.
\item $\NF_B(y)=\NF_B(z)$.
\end{enumerate}
Thus Definition~\ref{def:wall-graph} gives an exact finite classification
of all type-$B$ $\mgraph$-molecules in every rank and every fixed-$k$
block.
\end{theorem}

\begin{proof}
By Theorem~\ref{thm:all-edges}, every bidirected canonical edge is a strict
carrier edge.  By Theorem~\ref{thm:positive-local}, every positive-color
edge stays inside one occurrence.  By
Proposition~\ref{prop:tableau-wall-criterion}, every wall edge projects to
an edge of $\cO_{n,k}$.  Hence a molecule path projects to an occurrence--wall path
and preserves $\NF_B$.

Conversely, take a path in $\cO_{n,k}$ and choose a wall witness for each
of its edges.  Consecutive witnesses belonging to the same occurrence are
joined by the explicit moves $d_0^{\blacktriangledown},d_2,\ldots,d_{n-1}$ and
therefore, by Theorem~\ref{thm:positive-local}, by positive-color
bidirected carrier edges.  Insert these positive paths between successive
wall witnesses.  This lifts the occurrence--wall path to a molecule path
and attaches the prescribed endpoints.

Every $\tau_i$, including $\tau_0$, changes only signs or labels of
existing carriers.  It preserves the number $k$ of pair components, so
the whole lifted path remains in $\cX_{n,k}$.
\end{proof}

\begin{corollary}
There is no hidden canonical-basis computation in the final contraction:
once both tableaux of every candidate wall pair are known, molecules can
be computed using only the explicit tableau moves
$d_0^{\blacktriangledown},d_2,\ldots,d_{n-1}$ and the tableau descent
difference in \eqref{eq:wall-graph-tableau-edge}.  No inverse insertion
is used in the wall decision itself.
\end{corollary}

\begin{example}[A rank-six occurrence--wall graph]
\label{ex:rank-six-graph}
Take $n=6$ and $k=2$.  Four modified Beissinger occurrences have the
following representatives and cardinalities:
\begin{align*}
 \mathsf L&=\Occ(Z_{12}^+A_3^+A_4^-Z_{56}^+),&|\mathsf L|&=5,\\
 \mathsf D&=\Occ(A_1^-Z_{23}^+A_4^-Z_{56}^+),&|\mathsf D|&=16,\\
 \mathsf U&=\Occ(A_1^+Z_{23}^+A_4^-Z_{56}^+),&|\mathsf U|&=9,\\
 \mathsf R&=\Occ(Z_{12}^+A_3^-Z_{45}^+A_6^-),&|\mathsf R|&=5.
\end{align*}
The only occurrence--wall adjacency among these four vertices is
$\mathsf D--\mathsf U$.  It has five signed-matching witnesses.  Their
negative endpoints are
\begin{gather*}
 A_1^-Z_{24}^+Z_{36}^+A_5^-,\qquad
 A_1^-Z_{23}^+Z_{45}^+A_6^-,\qquad
 A_1^-Z_{23}^+Z_{46}^+A_5^-,\\
 A_1^-Z_{24}^+Z_{35}^+A_6^-,\qquad
 A_1^-Z_{23}^+A_4^-Z_{56}^+,
\end{gather*}
and the positive endpoint of each witness is obtained by replacing
$A_1^-$ by $A_1^+$.  For every such endpoint pair $(z^-,z^+)$,
\[
 1\in\Des\bigl(\Psi(z^+)\bigr)
       \setminus\Des\bigl(\Psi(z^-)\bigr),
\]
so all five adjacencies are detected directly by the tableau wall test.

\begin{center}
\begin{tikzpicture}[
  every node/.style={draw,rounded corners,align=center,
    minimum width=2.25cm,minimum height=1.05cm},
  x=1cm,y=1cm]
  \node (L) at (0,0) {$\mathsf L$\\[-1pt]\scriptsize $5$ vertices};
  \node (D) at (3.2,0) {$\mathsf D$\\[-1pt]\scriptsize $16$ vertices};
  \node (U) at (6.4,0) {$\mathsf U$\\[-1pt]\scriptsize $9$ vertices};
  \node (R) at (9.6,0) {$\mathsf R$\\[-1pt]\scriptsize $5$ vertices};
  \draw[very thick] (D)--node[draw=none,fill=white,yshift=18pt,
    inner sep=1pt,minimum width=0pt,minimum height=0pt]
    {\scriptsize five wall witnesses} (U);
\end{tikzpicture}

\smallskip
{\small Three connected components of $\cO_{6,2}$: the isolated
occurrences $\mathsf L$ and $\mathsf R$, and
$\mathsf D\cup\mathsf U$.}
\end{center}

The isolated occurrence $\mathsf L$ is already one entire molecule.  To
display a complete molecule rather than only its contracted occurrence,
number its five signed matchings as follows.  In the last column, the
positive component of $\Psi$ is displayed first and the negative-transpose
component second.
{\small
\[
\begin{array}{c@{\quad}l@{\qquad}l@{\qquad}c}
 &z&\Des_B^{\mathrm{tab}}(\Psi(z))&\Psi(z)\\ \hline
1&Z_{12}^+A_3^+A_4^-Z_{56}^+&\{0,2,4\}&
  \rowtableau{1&3\cr2&5\cr6\cr}\sqcup\rowtableau{4\cr}\\
2&Z_{12}^+A_3^+Z_{46}^+A_5^-&\{0,2,5\}&
  \rowtableau{1&3\cr2&4\cr6\cr}\sqcup\rowtableau{5\cr}\\
3&Z_{12}^+A_3^-A_4^+Z_{56}^+&\{0,3\}&
  \rowtableau{1&4\cr2&5\cr6\cr}\sqcup\rowtableau{3\cr}\\
4&Z_{13}^+A_2^+A_4^-Z_{56}^+&\{0,1,4\}&
  \rowtableau{1&2\cr3&5\cr6\cr}\sqcup\rowtableau{4\cr}\\
5&Z_{13}^+A_2^+Z_{46}^+A_5^-&\{0,1,3,5\}&
  \rowtableau{1&2\cr3&4\cr6\cr}\sqcup\rowtableau{5\cr}.
\end{array}
\]
}
Its complete bidirected graph is
\begin{center}
\begin{tikzpicture}[
  molvertex/.style={circle,draw,thick,minimum size=7.5mm,inner sep=0pt},
  mollabel/.style={draw=none,fill=white,inner sep=1pt,font=\scriptsize}]
  \node[molvertex] (m5) at (0,2.4) {$5$};
  \node[molvertex] (m2) at (-2.5,1.15) {$2$};
  \node[molvertex] (m4) at (2.5,1.15) {$4$};
  \node[molvertex] (m1) at (0,0) {$1$};
  \node[molvertex] (m3) at (0,-1.8) {$3$};
  \draw[thick] (m1)--node[mollabel,above left] {$s_4$} (m2);
  \draw[thick] (m1)--node[mollabel,above right] {$s_2$} (m4);
  \draw[thick] (m2)--node[mollabel,above right] {$s_2$} (m5);
  \draw[thick] (m4)--node[mollabel,above left] {$s_4$} (m5);
  \draw[thick] (m1)--node[mollabel,right] {$s_3$} (m3);
\end{tikzpicture}

{\small The complete five-vertex molecule $\mathsf L$; an edge labeled
$s_i$ is the bidirected carrier move of color $i$.}
\end{center}
The displayed five edges are all its bidirected edges: in particular,
there is no suppressed wall edge or additional positive-color edge.

Consequently these components expand to three molecules containing
$5$, $16+9=25$, and $5$ signed matchings, respectively.  The figure also
illustrates why the occurrence graph must remember the wall incidence:
the positive-color moves alone distinguish $\mathsf D$ and $\mathsf U$,
whereas the wall toggle joins them.
\end{example}

\begin{remark}[Exact algorithm versus closed statistic]
The invariant $\NF_B$ is an exact finite quotient presentation.  We have
not claimed a monotone rewriting system or a closed partition-valued
statistic analogous to ordinary Robinson--Schensted shape.  This
distinction prevents the connectivity step from being hidden inside the
word ``normal form.''
\end{remark}

The classification theorem is constructive.  To pass from the
occurrence--wall description to an explicit finite computation, apply the
following terminating procedure.

\begin{enumerate}[label=\textbf{Step \arabic*.},leftmargin=5.2em]
\item Enumerate the signed matchings in $\cX_{n,k}$.
\item For each $z$, compute $\PRBSp(z)$ by restricted row Beissinger
insertion and $\PCBSm(z)$ by restricted column Beissinger insertion, and
form
$\Psi(z)=\PRBSp(z)\sqcup\bigl(\PCBSm(z)\bigr)^{\mathsf t}$.
\item Form the components of $\cD_{n,k}$ using every nontrivial
$d_0^{\blacktriangledown},d_2,\ldots,d_{n-1}$.
\item Enumerate both endpoints of each candidate wall toggle and add an
edge of $\cO_{n,k}$ exactly when the tableau descent difference in
\eqref{eq:wall-graph-tableau-edge} is nonempty.
\item Return the connected components of $\cO_{n,k}$, with every
occurrence expanded back to its signed matchings.
\end{enumerate}

Termination is immediate because $\cX_{n,k}$ is finite.  Correctness is
Theorem~\ref{thm:classification}.

\section{The companion \texorpdfstring{$\ngraph$}{n}-classification}
\label{sec:n-case}

The $\ngraph$-case has the same positive-color combinatorics as the
$\mgraph$-case after all carrier signs are reversed, but its wall rule is
not Theorem~\ref{thm:wall} verbatim.  We first record the duality that
controls this distinction.  For $z\in\cX_{n,k}$, let $-z$ be obtained by
replacing every carrier $A_i^\epsilon$ and $Z_{ij}^\epsilon$ by
$A_i^{-\epsilon}$ and $Z_{ij}^{-\epsilon}$, respectively.

\begin{lemma}[Sign-opposition duality]
\label{lem:mn-duality}
For $x,y\in\cX_{n,k}$,
\begin{equation}
 x\longrightarrow y\text{ in the $\mgraph$-graph}
 \quad\Longleftrightarrow\quad
 -y\longrightarrow -x\text{ in the $\ngraph$-graph}.
 \label{eq:arrow-duality}
\end{equation}
Consequently,
\begin{equation}
 x\sim_{\mathrm{mol},\mgraph}y
 \quad\Longleftrightarrow\quad
 -x\sim_{\mathrm{mol},\ngraph}-y.
 \label{eq:molecule-duality}
\end{equation}
The descent labels satisfy \eqref{eq:mn-label-duality}.
\end{lemma}

\begin{proof}
Let $\rho_{\mgraph}$ and $\rho_{\ngraph}$ denote the two Gelfand
representations and let $\Theta$ be the Hecke involution
$\Theta(H_i)=-H_i^{-1}$.  The sign-opposition pairing of
\cite{MZGelfand} gives a signed permutation matrix $Q$, taking the
$z$-coordinate to the $-z$-coordinate, such that
\[
 \rho_{\ngraph}(-H_i^{-1})
   =Q\rho_{\mgraph}(H_i)^{\mathsf t}Q^{-1}
 \qquad(0\leq i<n).
\]
The quadratic relation gives
\[
                    -H_i^{-1}=-H_i+(v-v^{-1})I.
\]
The scalar summand has no off-diagonal entries.  Hence transposition and
sign reversal exchange the support of the coefficient from $x$ to $y$
with the support from $-y$ to $-x$, proving
\eqref{eq:arrow-duality}.  Applying the equivalence in both directions
to a bidirected pair proves \eqref{eq:molecule-duality}.  The diagonal
terms exchange ascent and descent labels, which is exactly
\eqref{eq:mn-label-duality}.
\end{proof}

The duality gives a direct insertion rule, rather than merely an abstract
bijection of molecules.

\begin{definition}[The split $\ngraph$-tableau]
\label{def:n-Psi}
Define
\begin{equation}
                         \Psi_{\ngraph}(z)=\Psi(-z).
\label{eq:n-Psi}
\end{equation}
Equivalently, apply restricted row Beissinger insertion to the negative
submatching of $z$, and apply restricted column Beissinger insertion
followed by transposition to the positive submatching.  In the ordered
row-reading word, the positive-transpose component precedes the
negative-row component.  The negative fixed points are inserted in
increasing order in the row construction; the positive fixed points are
inserted in decreasing order in the column construction before
transposition.
\end{definition}

Thus the $\ngraph$-construction literally exchanges the two sign roles in
Definition~\ref{def:Psi}.  Let $\cD^{\ngraph}_{n,k}$ be the graph on the
admissible tableaux $\Psi_{\ngraph}(z)$ whose edges are the nontrivial
moves $d_0^{\blacktriangledown},d_2,\ldots,d_{n-1}$, and put
\[
 \Occ_{\ngraph}(z)
   =[\Psi_{\ngraph}(z)]_{\cD^{\ngraph}_{n,k}}.
\]
Since negation commutes with every positive-color carrier conjugation,
Theorem~\ref{thm:positive-local} and Lemma~\ref{lem:mn-duality} show that
these are exactly the components obtained by deleting the wall color
from the $\ngraph$ bidirected graph.

\begin{proposition}[$\ngraph$ wall criterion]
\label{prop:n-wall}
The wall toggle $z\leftrightarrow\tau_0z$ is bidirected in the
$\ngraph$-graph if and only if one of the following mutually exclusive
conditions holds.
\begin{enumerate}[leftmargin=3.2em]
\item[\textup{(NF)}]
The component containing $1$ is a fixed point, $n\geq2$, and $2$ is not
a negative fixed point.
\item[\textup{(NP)}]
The component containing $1$ is $Z_{1r}^{\epsilon}$ and the component
containing $2$ is $Z_{2q}^{\delta}$ with $2<q<r$; the signs are arbitrary.
\end{enumerate}
In particular, there is no wall edge for $n=1$ or for
$Z_{12}^{\epsilon}$.
\end{proposition}

\begin{proof}
Negation commutes with $\tau_0$.  By Lemma~\ref{lem:mn-duality}, the
toggle is bidirected in the $\ngraph$-graph exactly when the toggle of
$-z$ is bidirected in the $\mgraph$-graph.  Condition \textup{(WP)} is
unchanged by reversing both arc signs.  In condition \textup{(WF)}, the
statement that $2$ is not positive fixed for $-z$ is precisely the
statement that $2$ is not negative fixed for $z$.  The result follows
from Theorem~\ref{thm:wall}.
\end{proof}

\begin{proposition}[Endpoint tableau form of the $\ngraph$-wall criterion]
\label{prop:n-tableau-wall-criterion}
For a split $\ngraph$-tableau $\Psi_{\ngraph}(u)$, set
\[
 \Des_{\ngraph}^{\mathrm{tab}}(\Psi_{\ngraph}(u))
   =[n-1]\setminus\Des(\Psi_{\ngraph}(u)).
\]
If $z^-,z^+$ are the negative and positive endpoints of a wall toggle,
respectively, then the endpoint tableau form of the $\ngraph$ criterion is
\begin{align}
 z^-\longleftrightarrow z^+
 &\quad\Longleftrightarrow\quad
 \Des_{\ngraph}^{\mathrm{tab}}(\Psi_{\ngraph}(z^+))
 \setminus
 \Des_{\ngraph}^{\mathrm{tab}}(\Psi_{\ngraph}(z^-))
 \ne\varnothing \notag\\
 &\quad\Longleftrightarrow\quad
 \Des(\Psi_{\ngraph}(z^-))
 \setminus\Des(\Psi_{\ngraph}(z^+))\ne\varnothing.
 \label{eq:n-tableau-wall-criterion}
\end{align}
As in the $\mgraph$ case, this criterion compares the two known endpoint
tableaux and does not construct one from the other.
\end{proposition}

\begin{proof}
Equations~\eqref{eq:mn-label-duality} and \eqref{eq:n-Psi}, together with
Lemma~\ref{lem:tableau-descents} identify
$\Des_{\ngraph}^{\mathrm{tab}}$ with the positive-node part of the
$\ngraph$ descent label.
The equivalence therefore follows from Proposition~\ref{prop:n-wall} and
the incomparability criterion of Theorem~\ref{thm:all-edges}.
\end{proof}

Let $\cO^{\ngraph}_{n,k}$ be the graph whose vertices are the
$\ngraph$-occurrences.  Join two occurrences when they contain the two
endpoints $z^-,z^+$ of a wall toggle satisfying
Proposition~\ref{prop:n-tableau-wall-criterion}.  Multiple endpoint
witnesses give one undirected edge.  As in Section~\ref{sec:classification}, fix a total
order on tableau pairs and let $\NF_B^{\ngraph}(z)$ denote the least
occurrence in the connected component of $\cO^{\ngraph}_{n,k}$ containing
$\Occ_{\ngraph}(z)$.

\begin{theorem}[Type-$B$ $\ngraph$-molecule classification]
\label{thm:n-classification}
For $y,z\in\cX_{n,k}$, the following are equivalent.
\begin{enumerate}[label=\textup{(\roman*)}]
\item The vertices $y,z$ belong to the same molecule of the type-$B$
Gelfand $\ngraph$-graph.
\item Their split $\ngraph$-occurrences lie in the same connected
component of $\cO^{\ngraph}_{n,k}$.
\item $\NF_B^{\ngraph}(y)=\NF_B^{\ngraph}(z)$.
\item The vertices $-y,-z$ belong to the same molecule of the type-$B$
Gelfand $\mgraph$-graph.
\end{enumerate}
Thus the $\ngraph$-molecules admit the same finite occurrence--wall
algorithm as in Section~\ref{sec:classification}, with
$\Psi_{\ngraph}$ and the reversed tableau descent difference
\eqref{eq:n-tableau-wall-criterion} in place of $\Psi$ and
\eqref{eq:tableau-wall-criterion}.
\end{theorem}

\begin{proof}
By \eqref{eq:n-Psi}, negation identifies
$\cD^{\ngraph}_{n,k}$ with $\cD_{n,k}$.  Proposition~\ref{prop:n-wall}
and Theorem~\ref{thm:wall} show that it also identifies
$\cO^{\ngraph}_{n,k}$ with $\cO_{n,k}$.  Theorem~\ref{thm:classification}
and Lemma~\ref{lem:mn-duality} now give all four equivalences.  The
algorithmic assertion follows because negation preserves the number of
pair components.
\end{proof}

\begin{example}[The rank-two polarity change]
In the $\ngraph$-graph one has the wall edge
\[
                 A_1^-A_2^+\longleftrightarrow A_1^+A_2^+,
\]
because $2$ is not a negative fixed point.  In contrast,
\[
                 A_1^-A_2^-\longleftrightarrow A_1^+A_2^-
\]
is not a $\ngraph$ wall edge.  The second toggle is the $\mgraph$ wall
edge from Theorem~\ref{thm:wall}; the first is its sign-opposite image.
This is the smallest example showing why the $\mgraph$ wall criterion
cannot be copied verbatim into the $\ngraph$ case.
\begin{center}
\begin{tikzpicture}[x=.8cm,y=.65cm,
  v/.style={circle,draw,fill=white,inner sep=0pt,minimum size=4.2mm},
  note/.style={font=\scriptsize,anchor=west}]
  \node[v] (a1) at (0,0) {\scriptsize$1$};
  \node[v] (a2) at (1,0) {\scriptsize$2$};
  \node[font=\scriptsize] at (0,.52) {$-$};
  \node[font=\scriptsize] at (1,.52) {$+$};
  \node at (1.75,0) {$\longleftrightarrow$};
  \node[v] (b1) at (2.5,0) {\scriptsize$1$};
  \node[v] (b2) at (3.5,0) {\scriptsize$2$};
  \node[font=\scriptsize] at (2.5,.52) {$+$};
  \node[font=\scriptsize] at (3.5,.52) {$+$};
  \node[note] at (4.15,0) {$\ngraph$ wall edge};

  \node[v] (c1) at (0,-1.25) {\scriptsize$1$};
  \node[v] (c2) at (1,-1.25) {\scriptsize$2$};
  \node[font=\scriptsize] at (0,-.73) {$-$};
  \node[font=\scriptsize] at (1,-.73) {$-$};
  \node at (1.75,-1.25) {$\nleftrightarrow$};
  \node[v] (d1) at (2.5,-1.25) {\scriptsize$1$};
  \node[v] (d2) at (3.5,-1.25) {\scriptsize$2$};
  \node[font=\scriptsize] at (2.5,-.73) {$+$};
  \node[font=\scriptsize] at (3.5,-.73) {$-$};
  \node[note] at (4.15,-1.25) {not an $\ngraph$ wall edge};
\end{tikzpicture}
\end{center}
\end{example}

\begin{example}[The rank-six $\ngraph$ occurrence--wall graph]
\label{ex:rank-six-n-graph}
Continue with $n=6$ and $k=2$.  Negating the four representatives in
Example~\ref{ex:rank-six-graph} gives the $\ngraph$-occurrences
\begin{align*}
 \mathsf L_{\ngraph}
   &=\Occ_{\ngraph}(Z_{12}^-A_3^-A_4^+Z_{56}^-),
 &|\mathsf L_{\ngraph}|&=5,\\
 \mathsf D_{\ngraph}
   &=\Occ_{\ngraph}(A_1^+Z_{23}^-A_4^+Z_{56}^-),
 &|\mathsf D_{\ngraph}|&=16,\\
 \mathsf U_{\ngraph}
   &=\Occ_{\ngraph}(A_1^-Z_{23}^-A_4^+Z_{56}^-),
 &|\mathsf U_{\ngraph}|&=9,\\
 \mathsf R_{\ngraph}
   &=\Occ_{\ngraph}(Z_{12}^-A_3^+Z_{45}^-A_6^+),
 &|\mathsf R_{\ngraph}|&=5.
\end{align*}
The only occurrence--wall adjacency among them is
$\mathsf D_{\ngraph}--\mathsf U_{\ngraph}$.  Its five positive endpoints
are
\begin{gather*}
 A_1^+Z_{24}^-Z_{36}^-A_5^+,
 \qquad A_1^+Z_{23}^-Z_{45}^-A_6^+,
 \qquad A_1^+Z_{23}^-Z_{46}^-A_5^+,\\
 A_1^+Z_{24}^-Z_{35}^-A_6^+,
 \qquad A_1^+Z_{23}^-A_4^+Z_{56}^-.
\end{gather*}
The other endpoint of each wall edge is obtained by changing $A_1^+$ to
$A_1^-$.  These are exactly the sign opposites of the five witnesses in
Example~\ref{ex:rank-six-graph}.  Each pair satisfies the reversed
tableau descent difference in \eqref{eq:n-tableau-wall-criterion}; hence
the five quotient adjacencies are obtained without reusing the
$\mgraph$ carrier test.

\begin{center}
\begin{tikzpicture}[
  every node/.style={draw,rounded corners,align=center,
    minimum width=2.25cm,minimum height=1.05cm},
  x=1cm,y=1cm]
  \node (L) at (0,0) {$\mathsf L_{\ngraph}$\\[-1pt]
    \scriptsize $5$ vertices};
  \node (D) at (3.2,0) {$\mathsf D_{\ngraph}$\\[-1pt]
    \scriptsize $16$ vertices};
  \node (U) at (6.4,0) {$\mathsf U_{\ngraph}$\\[-1pt]
    \scriptsize $9$ vertices};
  \node (R) at (9.6,0) {$\mathsf R_{\ngraph}$\\[-1pt]
    \scriptsize $5$ vertices};
  \draw[very thick] (D)--node[draw=none,fill=white,yshift=18pt,
    inner sep=1pt,minimum width=0pt,minimum height=0pt]
    {\scriptsize five $\ngraph$ wall witnesses} (U);
\end{tikzpicture}

\smallskip
{\small Three connected components of $\cO^{\ngraph}_{6,2}$: the
isolated occurrences $\mathsf L_{\ngraph}$ and
$\mathsf R_{\ngraph}$, and
$\mathsf D_{\ngraph}\cup\mathsf U_{\ngraph}$.}
\end{center}

The isolated occurrence $\mathsf L_{\ngraph}$ is the following complete
five-vertex $\ngraph$-molecule.  In the last column, the negative-row
component of $\Psi_{\ngraph}$ is displayed first and the
positive-transpose component second:
{\small
\[
\begin{array}{c@{\quad}l@{\qquad}l@{\qquad}c}
 &z&\Des_{\ngraph}^{\mathrm{tab}}(\Psi_{\ngraph}(z))
      &\Psi_{\ngraph}(z)\\ \hline
1&Z_{12}^-A_3^-A_4^+Z_{56}^-&\{1,3,5\}&
  \rowtableau{1&3\cr2&5\cr6\cr}\sqcup\rowtableau{4\cr}\\
2&Z_{12}^-A_3^-Z_{46}^-A_5^+&\{1,3,4\}&
  \rowtableau{1&3\cr2&4\cr6\cr}\sqcup\rowtableau{5\cr}\\
3&Z_{12}^-A_3^+A_4^-Z_{56}^-&\{1,2,4,5\}&
  \rowtableau{1&4\cr2&5\cr6\cr}\sqcup\rowtableau{3\cr}\\
4&Z_{13}^-A_2^-A_4^+Z_{56}^-&\{2,3,5\}&
  \rowtableau{1&2\cr3&5\cr6\cr}\sqcup\rowtableau{4\cr}\\
5&Z_{13}^-A_2^-Z_{46}^-A_5^+&\{2,4\}&
  \rowtableau{1&2\cr3&4\cr6\cr}\sqcup\rowtableau{5\cr}.
\end{array}
\]
}
Indeed, these augmented tableau label sets are the complements of the
corresponding $\mgraph$ tableau label sets in
Example~\ref{ex:rank-six-graph}, as required by
\eqref{eq:mn-label-duality}.  Its complete bidirected graph is
\begin{center}
\begin{tikzpicture}[
  molvertex/.style={circle,draw,thick,minimum size=7.5mm,inner sep=0pt},
  mollabel/.style={draw=none,fill=white,inner sep=1pt,font=\scriptsize}]
  \node[molvertex] (m5) at (0,2.4) {$5$};
  \node[molvertex] (m2) at (-2.5,1.15) {$2$};
  \node[molvertex] (m4) at (2.5,1.15) {$4$};
  \node[molvertex] (m1) at (0,0) {$1$};
  \node[molvertex] (m3) at (0,-1.8) {$3$};
  \draw[thick] (m1)--node[mollabel,above left] {$s_4$} (m2);
  \draw[thick] (m1)--node[mollabel,above right] {$s_2$} (m4);
  \draw[thick] (m2)--node[mollabel,above right] {$s_2$} (m5);
  \draw[thick] (m4)--node[mollabel,above left] {$s_4$} (m5);
  \draw[thick] (m1)--node[mollabel,right] {$s_3$} (m3);
\end{tikzpicture}

{\small The complete molecule $\mathsf L_{\ngraph}$; the graph is the
sign-opposite image of the complete molecule in
Example~\ref{ex:rank-six-graph}.}
\end{center}
Thus the four displayed occurrences expand to three $\ngraph$-molecules
of sizes $5$, $25$, and $5$, and the graph makes the sign-opposition
correspondence of Theorem~\ref{thm:n-classification} visible at the level
of every bidirected edge.
\end{example}

\section{The type-\texorpdfstring{$D$}{D} occurrence--fork classification}
\label{sec:type-d}

We now pass from $B_n$ to its even signed subgroup $D_n$.  Throughout this
section $n\geq3$, $0\leq k\leq\lfloor n/2\rfloor$, and all vertices lie
in the parity-compatible set $\cX^D_{n,k}$ of
\eqref{eq:type-D-carriers}.  The argument applies to the canonical
type-$D$ graphs in every rank.  The word ``Gelfand'' is used only when
$n$ is odd.

\subsection{Completion-filtered tail occurrences}

The split tableau $\Psi(z)$ of Definition~\ref{def:Psi} remains a
bijection from signed carriers to its image after restriction to
$\cX^D_{n,k}$.  It is useful as a coordinate system, but the type-$B$
positive-edge graph cannot simply be restricted to this subset.  Indeed,
the type-$D$ node $s_2$ is adjacent to both fork nodes, whereas the
type-$B$ wall belongs to a different Coxeter diagram.  We therefore retain
the fork bit from the full completion, but express every positive-node
label as a tableau descent.  Namely, define the augmented tableau descent
set
\begin{equation}
 \Des_D^{\mathrm{tab}}(\Psi(u))
   =\Des(\Psi(u))\cup
   \begin{cases}
    \{-1\},&s_{-1}\in D_D(u),\\
    \varnothing,&s_{-1}\notin D_D(u).
   \end{cases}
 \label{eq:D-augmented-tableau-descents}
\end{equation}
The fork bit on the right is computed by the four completion cases in
Section~\ref{subsec:type-D-prelim}.  Lemma~\ref{lem:tableau-descents}
and the equality of the positive-node type-$B$ and type-$D$ labels give
\[
                    \Des_D^{\mathrm{tab}}(\Psi(u))=D_D(u).
\]

\begin{definition}[Type-$D$ tail occurrence]
\label{def:D-tail-occurrence}
Let $\cD^D_{n,k}$ be the graph on the tableaux $\Psi(z)$ with
$z\in\cX^D_{n,k}$.  Join $\Psi(z)$ and $\Psi(z')$ when, for some
$1\leq i<n$,
\[
 z'=s_i z s_i\ne z,
\]
the conjugation is strict and the augmented tableau descent sets
$\Des_D^{\mathrm{tab}}(\Psi(z))$ and
$\Des_D^{\mathrm{tab}}(\Psi(z'))$ are incomparable.  Define
\[
 \Occ_D(z)=[\Psi(z)]_{\cD^D_{n,k}}.
\]
We call this component the \emph{completion-filtered type-$D$ tail
occurrence} of $z$.
\end{definition}

Every condition in this definition is finite and explicit.  The strict
and weak states at positive nodes are given by
\eqref{eq:descent-formula}; the extra label $s_{-1}$ is given in
Section~\ref{subsec:type-D-prelim}.  Thus no canonical-basis recursion is
part of the definition.

\begin{proposition}
\label{prop:D-tail}
The type-$D$ tail occurrences are exactly the connected components
obtained from the bidirected $\mgraph$-graph by deleting the extra fork
color $s_{-1}$.
\end{proposition}

\begin{proof}
Positive conjugation preserves the signs and the number of all carrier
components, hence preserves $\cX^D_{n,k}$.  By
Theorem~\ref{thm:all-edges}, a pair joined by a positive simple
conjugation is bidirected exactly when the conjugation is strict and its
two type-$D$ descent labels are incomparable.  By
\eqref{eq:D-augmented-tableau-descents}, these are precisely the tableau
edges in Definition~\ref{def:D-tail-occurrence}.
\end{proof}

\subsection{The fork criterion}

The key point is a general constraint on the labels of an edge in any
equal-parameter $W$-graph.

\begin{lemma}[Commuting-label lemma]
\label{lem:commuting-label}
Let $\Gamma=(V,\omega,\tau)$ be a $W$-graph over
$\mathbb Z[v,v^{-1}]$.  Suppose $x,y\in V$ satisfy
$\omega(y,x)\ne0$ and
\[
 r\in\tau(y)\setminus\tau(x),\qquad
 t\in\tau(x)\setminus\tau(y).
\]
Then the simple generators $r$ and $t$ do not commute.
\end{lemma}

\begin{proof}
Assume that $rt=tr$.  In the standard $W$-graph action,
$t\in\tau(x)$ gives $H_tC_x=-v^{-1}C_x$.  Since
$r\notin\tau(x)$ and $r\in\tau(y)$, the coefficient of $C_y$ in
$H_rC_x$ is $\omega(y,x)$.  Consequently, the coefficient of $C_y$ in
$H_rH_tC_x$ is
\[
                         -v^{-1}\omega(y,x).
\]
On the other hand, $t\notin\tau(y)$, so the diagonal coefficient of
$C_y$ in $H_tC_y$ is $v$.  No other term of $H_rC_x$ can contribute
$C_y$ after applying $H_t$: a non-diagonal occurrence of $C_y$ in
$H_tC_u$ would require $t\in\tau(y)$.  Hence the coefficient of $C_y$
in $H_tH_rC_x$ is
\[
                           v\omega(y,x).
\]
The relation $H_rH_t=H_tH_r$ would imply
$(v+v^{-1})\omega(y,x)=0$, contrary to $\omega(y,x)\ne0$.
\end{proof}

\begin{theorem}[Type-$D$ fork edge criterion]
\label{thm:D-fork}
Let $z,z'\in\cX^D_{n,k}$ satisfy $z'=\tau^D_{-1}z$.  Orient the pair
so that $\iota(z)$ is the lower endpoint and
$s_{-1}\iota(z)s_{-1}=\iota(z')$ is the upper endpoint of a strict cover.
Then
\begin{equation}
 z\longleftrightarrow z'
 \quad\Longleftrightarrow\quad
 s_2\in D_D(z)\setminus D_D(z').
 \label{eq:D-fork-criterion}
\end{equation}
\end{theorem}

\begin{proof}
At the lower endpoint, $s_{-1}$ is a strict ascent; at the upper endpoint
it is a strict descent.  Thus
\[
 s_{-1}\in D_D(z')\setminus D_D(z).
\]
If the edge is bidirected, Theorem~\ref{thm:all-edges} makes the two
descent labels incomparable.  Hence there is some
$t\in D_D(z)\setminus D_D(z')$.  Apply
Lemma~\ref{lem:commuting-label} with $r=s_{-1}$.  In the type-$D$ Coxeter
diagram the only simple generator not commuting with $s_{-1}$ is $s_2$;
therefore $t=s_2$.

Conversely, the displayed membership of $s_2$, together with the opposite
membership of $s_{-1}$, makes $D_D(z)$ and $D_D(z')$ incomparable.
Theorem~\ref{thm:all-edges} then makes the strict cover bidirected.
\end{proof}

\begin{proposition}[Endpoint tableau form of the fork criterion]
\label{prop:D-tableau-fork-criterion}
For the oriented endpoints $z<z'$ in Theorem~\ref{thm:D-fork}, one has
\begin{equation}
 z\longleftrightarrow z'
 \quad\Longleftrightarrow\quad
 2\in\Des(\Psi(z))\setminus\Des(\Psi(z')).
 \label{eq:D-tableau-fork-criterion}
\end{equation}
Thus, once the lower and upper endpoint tableaux are known, fork
bidirectedness is decided by the positions of the entries $2$ and $3$ in
the two split tableaux.  The orientation still comes from the strict
completion cover; the formula does not construct $\Psi(z')$ from
$\Psi(z)$.
\end{proposition}

\begin{proof}
By Lemma~\ref{lem:tableau-descents}, the condition
$s_2\in D_D(z)\setminus D_D(z')$ in Theorem~\ref{thm:D-fork} is
equivalent to
$2\in\Des(\Psi(z))\setminus\Des(\Psi(z'))$.
\end{proof}

\medskip
\noindent\emph{Illustrated fork example.}
In the type-$D_3$, $k=1$ block, orient the strict fork toggle as
\[
 z=Z_{23}^+A_1^+
 \mathrel{\mathop{\longrightarrow}^{s_{-1}}}
 z'=Z_{13}^-A_2^-.
\]
Its endpoint tableaux and ordinary tableau descents are
\[
\begin{array}{c@{\qquad}c@{\qquad}c}
u&\Psi(u)&\Des(\Psi(u))\\ \hline
z&\rowtableau{1\cr2\cr3\cr}\sqcup\varnothing&\{1,2\}\\[2pt]
z'&\varnothing\sqcup\rowtableau{1&3\cr2\cr}&\{1\}.
\end{array}
\]
The entry $3$ is below $2$ at the lower endpoint but not at the upper
endpoint, so $2\in\Des(\Psi(z))\setminus\Des(\Psi(z'))$.  Formula
\eqref{eq:D-tableau-fork-criterion} therefore turns the displayed strict
arrow into a bidirected fork edge.

Theorem~\ref{thm:D-fork} is a rank-two statement about $W$-graph labels,
but the endpoint tableaux themselves may still depend on the full
completion.  It does not assert that one endpoint tableau is determined
by the other, or that the fork edge is determined by the three carriers
through $1,2,3$.

\subsection{Classification and the companion graph}

\begin{definition}[Type-$D$ occurrence--fork graph]
\label{def:D-fork-graph}
Let $\cO^D_{n,k}$ have the type-$D$ tail occurrences as vertices.  Join
$L,L'$ when there are $z,z'\in\cX^D_{n,k}$ with
$\Psi(z)\in L$, $\Psi(z')\in L'$, and
$z'=\tau^D_{-1}z$, such that after orienting the strict cover as in
Theorem~\ref{thm:D-fork}, the endpoint tableaux satisfy
\begin{equation}
             2\in\Des(\Psi(z))\setminus\Des(\Psi(z')).
 \label{eq:D-fork-graph-tableau-edge}
\end{equation}
Multiple endpoint witnesses give one undirected edge.
\end{definition}

\begin{theorem}[Type-$D$ molecule classification]
\label{thm:D-classification}
For $x,y\in\cX^D_{n,k}$, the following are equivalent.
\begin{enumerate}[label=\textup{(\roman*)}]
\item $x$ and $y$ belong to the same molecule of
$\Gamma^{\mgraph}_{D,n,k}$.
\item The occurrences $\Occ_D(x)$ and $\Occ_D(y)$ lie in the same
connected component of $\cO^D_{n,k}$.
\item The vertices $\delta_n^D(x)$ and $\delta_n^D(y)$ belong to the
same molecule of $\Gamma^{\ngraph}_{D,n,k}$.
\end{enumerate}
When $n$ is odd these are classifications of the type-$D_n$ Gelfand
$\mgraph$- and $\ngraph$-molecules.
\end{theorem}

\begin{proof}
By Proposition~\ref{prop:D-tail}, every bidirected edge whose color is in
$\{s_1,\ldots,s_{n-1}\}$ stays inside one tail occurrence.  By
Proposition~\ref{prop:D-tableau-fork-criterion}, every
remaining bidirected edge projects to an edge of $\cO^D_{n,k}$.  Hence
every molecule path projects to an occurrence--fork path.

Conversely, choose a fork witness for each edge of an
occurrence--fork path.  Consecutive witnesses in one tail occurrence are
joined by a path of positive-color bidirected edges by the definition of
$\Occ_D$.  Inserting these tail paths between the chosen fork witnesses
lifts the quotient path to a bidirected path in the original graph.  This
proves the equivalence of \textup{(i)} and \textup{(ii)}.

Theorem~1.13 of \cite{MZGelfand} says that $\delta_n^D$ identifies the
$\mgraph$-graph with the dual of the parity-modified $\ngraph$-graph.
Duality reverses all arrows and complements the labels, so it preserves
the set of bidirected pairs.  The parity modification has the same
underlying directed graph as the original $\ngraph$-graph.  Therefore
$\delta_n^D$ induces a bijection of molecules, proving the equivalence
with \textup{(iii)}.
\end{proof}

\begin{example}[A complete type-$D_3$ $\mgraph$-molecule]
\label{ex:D3-m-molecule}
Consider the $k=1$ block.  It consists of the six signed carriers
\[
 Z_{12}^+A_3^+,\quad Z_{13}^+A_2^+,\quad Z_{23}^+A_1^+,\qquad
 Z_{12}^-A_3^-,\quad Z_{13}^-A_2^-,\quad Z_{23}^-A_1^-.
\]
Here and below we abbreviate $s_i$ by $i$ inside a descent-label set.
Applying the completion formulas and the four fork cases from
Section~\ref{subsec:type-D-prelim} gives
\[
\begin{array}{c@{\qquad}c@{\qquad}c}
z&\iota(z)(1),\ldots,\iota(z)(6)
  &\Des_D^{\mathrm{tab}}(\Psi(z))\\ \hline
Z_{12}^+A_3^+&(2,1,4,3,6,5)&\{-1,1\}\\
Z_{13}^+A_2^+&(3,4,1,2,6,5)&\{2\}\\
Z_{23}^+A_1^+&(4,3,2,1,6,5)&\{1,2\}\\
Z_{12}^-A_3^-&(-2,-1,-4,-3,6,5)&\{-1,1,2\}\\
Z_{13}^-A_2^-&(-3,-4,-1,-2,6,5)&\{-1,1\}\\
Z_{23}^-A_1^-&(-4,-3,-2,-1,6,5)&\{-1,2\}.
\end{array}
\tag{D3-$\mgraph$}
\]
The positive-color part has four tail occurrences:
\begin{align*}
 \mathsf P&=\{Z_{12}^+A_3^+,Z_{13}^+A_2^+\},
 &\mathsf K&=\{Z_{13}^-A_2^-,Z_{23}^-A_1^-\},\\
 \mathsf F&=\{Z_{23}^+A_1^+\},
 &\mathsf Q&=\{Z_{12}^-A_3^-\}.
\end{align*}
The two nontrivial tail edges have colors $s_2$ in $\mathsf P$ and
$s_1$ in $\mathsf K$.  There is exactly one fork adjacency, namely
\[
 Z_{13}^-A_2^-
       \mathrel{\mathop{\longleftrightarrow}^{s_{-1}}}
 Z_{23}^+A_1^+.
\]
Indeed, the two endpoint labels are $\{-1,1\}$ and $\{1,2\}$.
Equivalently, their ordinary tableau descent sets are $\{1\}$ and
$\{1,2\}$, as displayed in the illustrated endpoint test above.  After
the strict cover is oriented, tableau descent $2$ is lost at the lower
endpoint, exactly as in \eqref{eq:D-tableau-fork-criterion}.  Testing the
other strict conjugates with the same endpoint tableau test gives no
further bidirected edge.  Thus the complete occurrence--fork graph is
\begin{center}
\begin{tikzpicture}[
  every node/.style={draw,rounded corners,align=center,
    minimum width=2.15cm,minimum height=1.0cm},
  x=1cm,y=1cm]
  \node (P) at (0,0) {$\mathsf P$\\[-1pt]\scriptsize $2$ vertices};
  \node (K) at (3.1,0) {$\mathsf K$\\[-1pt]\scriptsize $2$ vertices};
  \node (F) at (6.2,0) {$\mathsf F$\\[-1pt]\scriptsize $1$ vertex};
  \node (Q) at (9.3,0) {$\mathsf Q$\\[-1pt]\scriptsize $1$ vertex};
  \draw[very thick] (K)--node[draw=none,fill=white,yshift=17pt,
    inner sep=1pt,minimum width=0pt,minimum height=0pt]
    {\scriptsize one $s_{-1}$ witness} (F);
\end{tikzpicture}

\smallskip
{\small The occurrence--fork graph in the type-$D_3$, $k=1$
$\mgraph$-block.}
\end{center}

The component $\mathsf K\cup\mathsf F$ expands to the following complete
three-vertex molecule.  The split tableau attached to each of its vertices
is
\[
\begin{array}{c@{\qquad}c}
z&\Psi(z)\\ \hline
Z_{13}^-A_2^-&
  \varnothing\sqcup\rowtableau{1&3\cr2\cr}\\
Z_{23}^+A_1^+&
  \rowtableau{1\cr2\cr3\cr}\sqcup\varnothing\\
Z_{23}^-A_1^-&
  \varnothing\sqcup\rowtableau{1&2\cr3\cr}.
\end{array}
\]
Here the two positions separated by $\sqcup$ retain the positive and
negative-transpose component order of Definition~\ref{def:Psi}.
\begin{center}
\begin{tikzpicture}[
  molvertex/.style={draw,rounded corners,thick,align=center,
    minimum width=3.15cm,minimum height=9mm},
  mollabel/.style={draw=none,fill=white,inner sep=1pt,font=\scriptsize}]
  \node[molvertex] (a) at (0,0)
    {$Z_{13}^-A_2^-$\\[-1pt]\scriptsize
      $\Des_D^{\rm tab}=\{-1,1\}$};
  \node[molvertex] (b) at (-4.2,-2.0)
    {$Z_{23}^+A_1^+$\\[-1pt]\scriptsize
      $\Des_D^{\rm tab}=\{1,2\}$};
  \node[molvertex] (c) at (4.2,-2.0)
    {$Z_{23}^-A_1^-$\\[-1pt]\scriptsize
      $\Des_D^{\rm tab}=\{-1,2\}$};
  \draw[very thick] (a)--node[mollabel,above left] {$s_{-1}$} (b);
  \draw[very thick] (a)--node[mollabel,above right] {$s_1$} (c);
\end{tikzpicture}

\smallskip
{\small A complete type-$D_3$ $\mgraph$-molecule.  The left edge is a
fork edge; the right edge lies inside the tail occurrence $\mathsf K$.}
\end{center}
Consequently the three molecules in this block have sizes $2$, $3$, and
$1$.  This example is the smallest one in which the occurrence--fork
contraction is nontrivial.
\end{example}

\begin{example}[The dual type-$D_3$ $\ngraph$-molecule]
\label{ex:D3-n-molecule}
In the same $k=1$ block, the parity rule in
\cite[(3.26)]{MZGelfand} gives $\delta_3^D(z)=-z$; no diagram
automorphism is needed.  Hence the four $\ngraph$ tail occurrences are
\begin{align*}
 \mathsf P_{\ngraph}&=-\mathsf P
   =\{Z_{12}^-A_3^-,Z_{13}^-A_2^-\},
 &\mathsf K_{\ngraph}&=-\mathsf K
   =\{Z_{13}^+A_2^+,Z_{23}^+A_1^+\},\\
 \mathsf F_{\ngraph}&=-\mathsf F=\{Z_{23}^-A_1^-\},
 &\mathsf Q_{\ngraph}&=-\mathsf Q=\{Z_{12}^+A_3^+\}.
\end{align*}
The only fork adjacency is
$\mathsf K_{\ngraph}--\mathsf F_{\ngraph}$.  Thus its contracted graph
is
\begin{center}
\begin{tikzpicture}[
  every node/.style={draw,rounded corners,align=center,
    minimum width=2.15cm,minimum height=1.0cm},
  x=1cm,y=1cm]
  \node (P) at (0,0) {$\mathsf P_{\ngraph}$\\[-1pt]
    \scriptsize $2$ vertices};
  \node (K) at (3.1,0) {$\mathsf K_{\ngraph}$\\[-1pt]
    \scriptsize $2$ vertices};
  \node (F) at (6.2,0) {$\mathsf F_{\ngraph}$\\[-1pt]
    \scriptsize $1$ vertex};
  \node (Q) at (9.3,0) {$\mathsf Q_{\ngraph}$\\[-1pt]
    \scriptsize $1$ vertex};
  \draw[very thick] (K)--node[draw=none,fill=white,yshift=17pt,
    inner sep=1pt,minimum width=0pt,minimum height=0pt]
    {\scriptsize one $s_{-1}$ witness} (F);
\end{tikzpicture}

\smallskip
{\small The occurrence--fork graph in the type-$D_3$, $k=1$
$\ngraph$-block.}
\end{center}

The descent labels of the three vertices in the nontrivial component are
the complements of the labels of their $\mgraph$ preimages in
Example~\ref{ex:D3-m-molecule}.  Their split $\ngraph$-tableaux are
\[
\begin{array}{c@{\qquad}c}
z&\Psi_{\ngraph}(z)\\ \hline
Z_{13}^+A_2^+&
  \varnothing\sqcup\rowtableau{1&3\cr2\cr}\\
Z_{23}^-A_1^-&
  \rowtableau{1\cr2\cr3\cr}\sqcup\varnothing\\
Z_{23}^+A_1^+&
  \varnothing\sqcup\rowtableau{1&2\cr3\cr}.
\end{array}
\]
Its full expansion is therefore
\begin{center}
\begin{tikzpicture}[
  molvertex/.style={draw,rounded corners,thick,align=center,
    minimum width=3.15cm,minimum height=9mm},
  mollabel/.style={draw=none,fill=white,inner sep=1pt,font=\scriptsize}]
  \node[molvertex] (a) at (0,0)
    {$Z_{13}^+A_2^+$\\[-1pt]\scriptsize
      $D_D^{\ngraph}=\{2\}$};
  \node[molvertex] (b) at (-4.2,-2.0)
    {$Z_{23}^-A_1^-$\\[-1pt]\scriptsize
      $D_D^{\ngraph}=\{-1\}$};
  \node[molvertex] (c) at (4.2,-2.0)
    {$Z_{23}^+A_1^+$\\[-1pt]\scriptsize
      $D_D^{\ngraph}=\{1\}$};
  \draw[very thick] (a)--node[mollabel,above left] {$s_{-1}$} (b);
  \draw[very thick] (a)--node[mollabel,above right] {$s_1$} (c);
\end{tikzpicture}

\smallskip
{\small A complete type-$D_3$ $\ngraph$-molecule.}
\end{center}
Every edge is the sign-opposite, arrow-reversed image of the corresponding
edge in Example~\ref{ex:D3-m-molecule}.  Since the edges are bidirected,
arrow reversal does not change the displayed undirected graph.  The
$\ngraph$ block therefore again has molecules of sizes $2$, $3$, and
$1$, as required by Theorem~\ref{thm:D-classification}.
\end{example}

\begin{remark}[Scope of the type-$D$ result]
Theorem~\ref{thm:D-classification} is a complete finite classification:
it uses only split insertion, the fork-orientation bit from the explicit
completion formulas, endpoint tableau descent comparisons, and graph
connectivity.  It does not use a cell classification.  The stronger
assertion that every fork decision is a
bounded three-carrier rule would require an additional locality theorem
and is not claimed here.  In particular, the external partner assigned to
a positive fixed point is shifted by the total number of negative fixed
points, while negative fixed points are numbered in decreasing order.
Distant fixed carriers can therefore change the exact values of
$\iota(z)(1),\iota(z)(2),\iota(z)(3)$.
\end{remark}


\begin{thebibliography}{99}

\bibitem{Assaf}
S.~H. Assaf,
\newblock Dual equivalence graphs I: a new paradigm for Schur positivity,
\newblock \emph{Forum Math. Sigma} \textbf{3} (2015), e12.

\bibitem{Beissinger}
J.~S. Beissinger,
\newblock Similar constructions for Young tableaux and involutions, and
their application to shiftable tableaux,
\newblock \emph{Discrete Math.} \textbf{67} (1987), 149--163.

\bibitem{BonnafeIancu}
C.~Bonnaf\'e and L.~Iancu,
\newblock Left cells in type $B_n$ with unequal parameters,
\newblock \emph{Represent. Theory} \textbf{7} (2003), 587--609.

\bibitem{Garfinkle}
D.~Garfinkle,
\newblock On the classification of primitive ideals for complex
classical Lie algebras, I,
\newblock \emph{Compositio Math.} \textbf{75} (1990), 135--169.

\bibitem{GarfinkleVogan}
D.~Garfinkle and D.~A. Vogan, Jr.,
\newblock On the structure of Kazhdan--Lusztig cells for branched Dynkin
diagrams,
\newblock \emph{J. Algebra} \textbf{153} (1992), 91--120.

\bibitem{Haiman}
M.~D. Haiman,
\newblock Dual equivalence with applications, including a conjecture of
Proctor,
\newblock \emph{Discrete Math.} \textbf{99} (1992), 79--113.

\bibitem{KazhdanLusztig}
D.~Kazhdan and G.~Lusztig,
\newblock Representations of Coxeter groups and Hecke algebras,
\newblock \emph{Invent. Math.} \textbf{53} (1979), 165--184.

\bibitem{LusztigUnequal}
G.~Lusztig,
\newblock \emph{Hecke Algebras with Unequal Parameters},
\newblock CRM Monograph Series, vol.~18, American Mathematical Society,
Providence, RI, 2003.

\bibitem{MarbergBar}
E.~Marberg,
\newblock Bar operators for quasiparabolic conjugacy classes in a Coxeter
group,
\newblock \emph{J. Algebra} \textbf{453} (2016), 325--363.

\bibitem{MZGelfand}
E.~Marberg and Y.~Zhang,
\newblock Gelfand $W$-graphs for classical Weyl groups,
\newblock \emph{J. Algebra} \textbf{609} (2022), 292--336.

\bibitem{MZPerfect}
E.~Marberg and Y.~Zhang,
\newblock Perfect models for finite Coxeter groups,
\newblock \emph{J. Pure Appl. Algebra} \textbf{227} (2023), 107303.

\bibitem{MZInsertion}
E.~Marberg and Y.~Zhang,
\newblock Insertion algorithms for Gelfand $S_n$-graphs,
\newblock \emph{Ann. Comb.} \textbf{28} (2024), 1199--1242.

\bibitem{McGovernDomino}
W.~M. McGovern,
\newblock Left cells and domino tableaux in classical Weyl groups,
\newblock \emph{Compositio Math.} \textbf{101} (1996), 77--98.

\bibitem{Pietraho}
T.~Pietraho,
\newblock Equivalence classes in the Weyl groups of type $B_n$,
\newblock \emph{J. Algebraic Combin.} \textbf{27} (2008), 247--262.

\bibitem{Stembridge}
J.~R. Stembridge,
\newblock Admissible $W$-graphs,
\newblock \emph{Represent. Theory} \textbf{12} (2008), 346--368.

\end{thebibliography}
\end{document}